\documentclass[dvips]{siamltex}

\usepackage{amsfonts}
\usepackage{amsmath}
\usepackage{amssymb}
\usepackage{mathrsfs}
\usepackage{graphicx}
\usepackage{psfrag}

\newtheorem{remark}{Remark}

\newcommand{\VVV}{{\mathcal V }}
\newcommand{\CCC}{{\mathcal C }}

\newcommand{\RR}{{\mathbb R}}
\newcommand{\NN}{{\mathbb N}}

\newcommand{\SSS}{{\mathbb  S}}

\newcommand{\CC}{{\mathbb C}}

\newcommand{\her}[2]{\left\langle#1~|~#2\right\rangle}

\title{Approximate stabilization of a quantum particle in a
1D infinite square potential well \thanks{This work was supported in
part by the "Agence Nationale de la Recherche" (ANR), Projet Blanc
CQUID number 06-3-13957. \newline A short and preliminary version of
this paper has been submitted to the IFAC world congress 2008. }}

\author{
    Karine Beauchard\thanks{CNRS, CMLA,  ENS Cachan, Avenue du pr\'{e}sidant Wilson, 94230
    Cachan, France ({\tt karine.beauchard@cmla.ens-cachan.fr}).}
    \and Mazyar Mirrahimi\thanks{INRIA Rocquencourt, Domaine de
    Voluceau, Rocquencourt B.P. 105, 78153 Le Chesnay Cedex, France,
    ({\tt mazyar.mirrahimi@inria.fr}).} }

\begin{document}
\bibliographystyle{plain}
\maketitle

\begin{abstract}
We consider a non relativistic charged particle in a 1D infinite
square potential well. This quantum system is subjected to a
control, which is a uniform (in space) time depending electric
field. It is represented by a complex probability amplitude solution
of a Schr\"odinger equation on a 1D bounded domain, with Dirichlet
boundary conditions. We prove the almost global approximate
stabilization of the eigenstates by explicit feedback laws.
\end{abstract}

\begin{keywords}
    control of partial differential equations, bilinear Schr\"odinger equation,
    quantum systems, Lyapunov stabilization
\end{keywords}

\begin{AMS}
    93c20, 
    35Q40, 
    93D15 
\end{AMS}


\section{Introduction}\label{sec:intro}
\subsection{Main result}
As in ~\cite{rouchon-ri02bis,beauchard-04,beauchard-coron-05}, we
consider a non-relativist charged particle in a one dimensional
space, with a potential $V(x)$, in a uniform electric field
$t\mapsto u(t)\in\RR$. Under the dipole moment approximation
assumption, and after appropriate changes of scales, the evolution
of the particle's wave function is given by the following
Schr\"{o}dinger equation
$$
i\frac{\partial \Psi}{\partial
t}(t,x)=-\frac{\partial^2\Psi}{\partial
x^2}(t,x)+(V(x)-u(t)x)\Psi(t,x).
$$
We study the case of an infinite square potential well: $V(x)=0$ for
$x\in I:=(-1/2,1/2)$ and $V(x)=+\infty$ for $x$ outside $I$.
Therefore our system is
\begin{alignat}{3}
i\frac{\partial \Psi}{\partial t}(t,x) &=
-\frac{1}{2}\frac{\partial^2\Psi}{\partial
x^2}(t,x)-u(t)x\Psi(t,x), \qquad x\in I \label{eq:main}\\
\Psi(0,x)&=\Psi_0(x), \label{eq:initial}\\
\Psi(t,\pm1/2)&=0. \label{eq:bords}
\end{alignat}
It is a nonlinear control system, denoted by $(\Sigma)$, in which
\begin{itemize}
\item the control is the electric field
$u : \RR_{+} \rightarrow \RR$,
\item the state is the wave function
$\Psi: \RR_{+} \times I \rightarrow \CC$ with
$\Psi(t) \in \SSS$ for every $t \geqslant 0$,
\end{itemize}
where $\SSS:=\{\varphi\in L^2(I;\CC) ; \|\varphi\|_{L^2}=1\}$.

Let us introduce the operator $A$ defined by
$$\begin{array}{ll}
D(A):=(H^{2} \cap H^{1}_{0})(I,\CC), &
A \varphi := - \frac{1}{2} \frac{d^{2}\varphi}{dx^{2}},
\end{array}$$
and for $s \in \RR$  the spaces
$$H^{s}_{(0)}(I,\CC):=D(A^{s/2}).$$
The following Proposition recalls classical existence and uniqueness
results for the solutions of
(\ref{eq:main})-(\ref{eq:initial})-(\ref{eq:bords}). For sake of
completeness, a proof of this proposition is given in the Appendix.

\begin{proposition} \label{Existence-classic}
Let $\Psi_{0} \in \SSS$, $T>0$ and $u \in C^{0}([0,T],\RR)$.
There exists a unique weak solution of
(\ref{eq:main})-(\ref{eq:initial})-(\ref{eq:bords}),
i.e. a function
$\Psi \in C^{0}([0,T],\SSS) \cap C^{1}([0,T],H^{-2}_{(0)}(I,\mathbb{C}))$
such that
\begin{equation} \label{weak}
\Psi(t)=e^{-iAt}\Psi_{0} +i\int_{0}^{t}
e^{-iA(t-s)} u(s) x \Psi(s) ds
\text{ in } L^{2}(I,\CC) \text{ for every } t \in [0,T].
\end{equation}
and then (\ref{eq:main}) holds in $H^{-2}_{(0)}(I,\CC)$
for every $t \in [0,T]$.

If, moreover, $\Psi_{0} \in (H^{2}\cap H^{1}_{0})(I,\CC)$, then
$\Psi$ is a strong solution i.e.
$\Psi \in C^{0}([0,T],(H^{2}\cap H^{1}_{0})(I,\CC)) \cap
C^{1}([0,T],L^{2}(I,\CC))$,
the equality (\ref{eq:main}) holds in $L^{2}(I,\CC)$, for every $t \in [0,T]$,
the equality (\ref{eq:initial}) holds in $H^{2} \cap H^{1}_{0}(I,\CC)$ and
the equality (\ref{eq:bords}) holds for every $t \in [0,T]$.

The weak (resp. strong) solutions are continuous with respect to initial
conditions for the $C^{0}([0,T],L^{2})$-topology
(resp. for the $C^{0}([0,T],H^{2} \cap H^{1}_{0})$-topology.)
\end{proposition}

The symbol $\langle . , . \rangle$ denotes the usual Hermitian
product of $L^{2}(I,\CC)$ i.e.
$$\langle \varphi , \xi \rangle :=  \int_{I} \varphi(x) \overline{\xi(x)} dx.$$
For $\sigma \in \mathbb{R}$, we introduce
the operator $A_{\sigma}$ defined by
$$\begin{array}{ll}
D(A_{\sigma}):=(H^{2}\cap H^{1}_{0})(I;\CC),
&
A_{\sigma} \varphi := - \frac{1}{2}
\frac{\partial^{2} \varphi}{dx^{2}} - \sigma x \varphi.
\end{array}$$
It is well known that there exists an orthonormal basis
$(\phi_{k,\sigma})_{k \in \NN^{*}}$ of $L^{2}(I,\CC)$ of eigenvectors
of $A_{\sigma}$ :
$$\begin{array}{ll}
\phi_{k,\sigma} \in H^{2} \cap H^{1}_{0}(I,\CC),
&
A_{\sigma}\phi_{k,\sigma}=\lambda_{k,\sigma} \phi_{k,\sigma}
\end{array}$$
where $(\lambda_{k,\sigma})_{k \in \mathbb{N}^{*}}$ is a
non-decreasing sequence of real numbers. For $s>0$ and $\sigma \in
\RR$, we define
$$H^{s}_{(\sigma)}(I,\CC):=D(A_{\sigma}^{s/2}),$$
equipped with the norm
$$\|\varphi\|_{H^{s}_{(\sigma)}} := \left( \sum_{k=1}^{\infty}
\lambda_{k,\sigma}^{s} |\langle \varphi , \phi_{k,\sigma} \rangle |^{2} \right)^{1/2}.$$
For $k \in \NN^{*}$ and $\sigma \in \RR$, we define
$$
\mathcal{C}_{k,\sigma}:=\{ \phi_{k,\sigma} e^{i\theta} ;
\theta \in [0,2\pi) \}.
$$
In order to simplify the notations, we will write $\phi_{k}$,
$\lambda_{k}$, $\mathcal{C}_{k}$ instead of $\phi_{k,0}$,
$\lambda_{k,0}$, $\mathcal{C}_{k,0}$. We have
\begin{equation}\label{eq:spectre}
\lambda_k=\frac{k^2\pi^2}{2},\quad \phi_k=
\left\{
\begin{array}{l}
\sqrt{2}\cos(k\pi x), \text{ when } k \text{ is odd,}
\\
\sqrt{2}\sin(k\pi x), \text{ when } k \text{ is even.}
\end{array}
\right.
\end{equation}
The goal of this paper is the study of the stabilization of the
system $(\Sigma)$ around the eigenstates $\phi_{k,\sigma}$. More
precisely, for $k \in \NN^{*}$ and $\sigma \in \RR$ small, we state
feedback laws $u=u_{k,\sigma}(\Psi)$ for which the solution of
(\ref{eq:main})-(\ref{eq:initial})-(\ref{eq:bords}) with
$u(t)=u_{k,\sigma}(\Psi(t))$ is such that
$$\limsup\limits_{t \rightarrow + \infty}
\text{dist}_{L^{2}(I,\CC)}(\Psi(t),\mathcal{C}_{k,\sigma})$$ is
arbitrarily small. We consider the convergence toward the circle
$\mathcal{C}_{k,\sigma}$ because the wave function $\Psi$ is defined
up to a phase factor. For simplicity sakes, we will only work with
the ground state $\phi_{1,\sigma}$. However, the whole arguments
remain valid for the general case.

Note that, even though the feedback stabilization of a quantum
system necessitates more complicated models taking into account the
measurement backaction on the system (see
e.g.~\cite{haroche-CDF,vanhandel-et-al-05,mirrahimi-vanhandel-05}),
the kind of strategy considered in this paper can be helpful for the
open-loop control of closed quantum systems. Indeed, one can apply
the stabilization techniques for the Schr\"odinger equation in
simulation and retrieve the control signal that will be then applied
in open-loop on the real physical system. As it will be detailed
below, in the bibliographic overview, such kind of strategy has been
widely used in the context of finite dimensional quantum systems.

The main result of this article is the following one.
\begin{theorem} \label{Stab-Hs-sigma}
Let $\Gamma>0$, $s>0$, $\epsilon >0$, $\gamma \in (0,1)$.
There exists $\sigma^{**}=\sigma^{**}(\Gamma,s)>0$ such that,
for every $\sigma \in (-\sigma^{**},\sigma^{**})$,
there exists a feedback law
$v_{\sigma,\Gamma,s,\epsilon,\gamma}(\Psi)$
such that, for every
$\Psi_{0} \in  \SSS \cap (H^{2} \cap H^{1}_{0} \cap H^{s}_{(\sigma)})(I,\CC)$
with
$$\|\Psi_{0}  \|_{H^{s}_{(\sigma)}} \leqslant \Gamma
\text{  and  }
|\langle \Psi_{0} , \phi_{1,\sigma} \rangle | > \gamma,$$
the Cauchy problem (\ref{eq:main})-(\ref{eq:initial})-(\ref{eq:bords})
with $u(t)=\sigma + v_{\sigma,\Gamma,s,\epsilon,\gamma}(\Psi)$
has a unique strong solution, moreover, this solution satisfies
$$
\limsup\limits_{t \rightarrow + \infty}
\text{dist}_{L^{2}} (\Psi(t) , \mathcal{C}_{1,\sigma}) \leqslant \epsilon.
$$
\end{theorem}

For $\sigma \neq 0$, the feedback law will be given explicitly. For
$\sigma = 0$, the feedback law will be given by an implicit formula.
The Theorem~\ref{Stab-Hs-sigma} provides \textit{almost global
approximate stabilization}. Indeed, any initial condition $\Psi_{0}
\in \SSS$ such that $\Psi_{0} \in H^{s}(I,\CC)$ for some $s>0$ and
$\langle \Psi_{0} , \phi_{1,\sigma} \rangle \neq 0$ can be moved
approximately to the circle $\mathcal{C}_{1,\sigma}$, thanks to an
appropriate feedback law. We will see that the assumption
\textquotedblleft$\Psi_{0} \in H^{s}(I,\CC)$, for some
$s>0$\textquotedblright is not necessary for doing that. In fact,
even for a $\Psi_{0}$ only belonging to $\SSS$, we can find the
appropriate feedback law as a function of the initial state
$\Psi_{0}$.

Notice that, physically, the assumption $\langle \Psi_{0} ,
\phi_{1,\sigma} \rangle \neq 0$ is not really restrictive. Indeed,
if $\langle \Psi_{0} , \phi_{1,\sigma} \rangle = 0$, a control field
in resonance with the natural frequencies of the system (the
difference between the eigenvalues corresponding to an eigenstate
whose population in the initial state is non-zero and the ground
state) will, instantaneously, ensure a non-zero population of the
ground state in the wavefunction. Then, one can just apply the
feedback law of the Theorem~\ref{Stab-Hs-sigma}.

\subsection{A brief bibliography}
The controllability of a finite dimensional quantum system,
$\iota\frac{d}{dt}\Psi=(H_0+u(t)~H_1)\Psi$ where $\Psi\in \CC^N$ and
$H_0$ and $H_1$ are $N\times N$ Hermitian matrices with coefficients
in $\CC$ has been very well
explored~\cite{sussmann-jurdjevic-72,ramakrishna-et-al-95,
albertini-et-al-03,altafini-JMP-02,turinici-rabitz-03}.
However, this does not guarantee the simplicity of the trajectory
generation. Very often the chemists formulate the task of the
open-loop control as a cost functional to be minimized. Optimal
control techniques (see e.g.,~\cite{hbref17}) and iterative
stochastic techniques (e.g, genetic
algorithms~\cite{turinici-control1}) are then two classes of
approaches which are most commonly used for this task.

When some non-degeneracy assumptions concerning the linearized
system are satisfied,~\cite{mirrahimi-et-al2-04} provides another
method based on Lyapunov techniques for generating trajectories. The
relevance of such a method for the control of chemical models has
been studied in~\cite{mirrahimi-et-al04}. As mentioned above, the
closed-loop system is simulated and the retrieved control signal is
applied in open-loop. Such kind of strategy has already been applied
widely in this framework~\cite{rabitz-tracking-95,sugawara-JCP03}.

The situation is much more difficult when we consider an infinite
dimensional configuration and less results are available. However,
the controllability of the system
(\ref{eq:main})-(\ref{eq:initial})-(\ref{eq:bords}) is now well
understood. In ~\cite{turinici-cdc00}, the author states some
non-controllability results for general Schr\"odinger systems. These
results apply in particular to the system
(\ref{eq:main})-(\ref{eq:initial})-(\ref{eq:bords}). However, this
negative result is due to the choice of the functional space that
does not allow the controllability. Indeed, if we consider different
functional spaces, one can get positive controllability results. In
~\cite{beauchard-04}, the local controllability of the system
(\ref{eq:main})-(\ref{eq:initial})-(\ref{eq:bords}) around the
ground state $\phi_{1,\sigma}$, for $\sigma$ small is proved. The
case $\sigma \neq 0$ is easier because the linearized system around
$\phi_{1,\sigma}$ for $\sigma \neq 0$ small is controllable; this
case is treated with the moment theory and a Nash-Moser implicit
functions theorem. As it has been discussed
in~\cite{rouchon-ri02bis}, the case of $\sigma=0$ is degenerate: the
linearized system around $\phi_1$ is not controllable. Therefore, in
this case, one needs to apply other tools, namely the return method
(introduced in ~\cite{coron-mcss-92}) and the quantum adiabatic
theory~\cite{avron-99}. In ~\cite{beauchard-coron-05}, the
steady-state controllability of this nonlinear system is proved
(i.e. the particle can be moved in finite time from an eigenstate
$\phi_{k}$ to another one $\phi_{j}$). The proof relies on many
local controllability results (proved with the previous strategies)
together with a compactness argument.

Concerning the trajectory generation problem for infinite
dimensional systems still much less results are available. The very
few existing literature is mostly based on the use of the optimal
control techniques~\cite{baudouin-et-al-05,baudouin-salomon-06}. The
simplicity of the feedback law found by the Lyapunov techniques
in~\cite{mirrahimi-et-al2-04,beauchard-et-al04} suggests the use of
the same approach for infinite dimensional configurations. However,
an extension of the convergence analysis to the PDE configuration is
not at all a trivial problem. Indeed, it requires the
pre-compactness of the closed-loop trajectories, a property that is
difficult to prove in infinite dimension. This strategy is used, for
exemple in \cite{JMC-BAN}.

In~\cite{mirrahimi-cdc06}, one of the authors proposes a
Lyapunov-based method to approximately stabilize a particle in a 3D
finite potential well under some restrictive assumptions. The author
assumes that the system is initialized in the finite dimensional
discrete part of the spectrum. Then, the idea consists in proposing
a Lyapunov function which encodes both the distance with respect to
the target state and the necessity of remaining in the discrete part
of the spectrum. In this way, he prevents the possibility of the
``mass lost phenomenon'' at infinity. Finally, applying some
dispersive estimates of Strichartz type, he ensures the approximate
stabilization of an arbitrary eigenstate in the discrete part of the
spectrum.

Finally, let us mention other strategies for proving the
stabilization of control systems. One can try to build a feedback
law for which one has a strict Lyapunov function. This strategy is
used, for example, for hyperbolic systems of conservation laws in
\cite{JMC-BAN-Lyapunov-strict}, for the 2-D incompressible Euler
equation in a simply connected domain in \cite{JMC-Stab-Euler2D},
see also \cite{Glass-Stab-Euler2D} for the multiconnected case. For
systems having a non controllable linearized system around the
equilibrium considered, the return method often provides good
results, see for example \cite{coron-mcss-92} for controllable
systems without drift and \cite{Glass-Camassa}) for Camassa-Holm
equation. In the end, we refer to \cite{JMC-book} for a pedagogical
presentation of strategies for the proof of stabilization of PDE
control systems.
\\

In this paper, we study the stabilization of the ground state
$\phi_{1,\sigma}$ for $\sigma$ in a neighborhood of $0$. Adapting
the techniques proposed in~\cite{mirrahimi-cdc06}, we ensure the
approximate stabilization of the system around $\phi_{1,\sigma}$.
Note that, the whole arguments hold if we replace the target state
by any eigenstate $\phi_{k,\sigma}$ of the system.

\subsection{Heuristic of the proof}\label{ssec:heur}
While trying to stabilize the ground state $\phi_{1,\sigma}$, a
first approach would be to consider the simple Lyapunov function
$$
\widetilde\VVV(\Psi)=1-|\her{\Psi}{\phi_{1,\sigma}}|^2.
$$
Just as in the finite the dimensional case~\cite{beauchard-et-al04},
the feedback law
$$
\tilde
u(\Psi)=\Im(\her{x\Psi}{\phi_{1,\sigma}}\her{\phi_{1,\sigma}}{\Psi})
$$
where $\Im$ denotes the imaginary part of a complex, ensures the
decrease of the Lyapunov function. However, trying to adapt the
convergence analysis, based on the use of the LaSalle invariance
principle, the pre-compactness of the trajectories in $L^2$
constitutes a major obstacle. Note that, in order to be able to
apply the LaSalle principle for an infinite dimensional system, one
certainly needs to prove such a pre-compactness result. In the
particular case of the infinite potential well, it even seems that,
one can not hope such a result. Indeed, phenomenons such as the
$L^2$-mass lost in the high energy levels do not allow this property
to hold true.

Similarly to~\cite{mirrahimi-cdc06}, the approach of this paper is
to avoid the population to go through the very high energy levels,
while trying to stabilize the system around $\phi_{1,\sigma}$.

As in Theorem~\ref{Stab-Hs-sigma}, let us consider $\Gamma>0$, $s
>0$, $\epsilon >0$,$\gamma >0$, $\sigma \in \mathbb{R}$. First, we
consider the case , $\sigma \neq 0$. Let $\Psi_0 \in
H^{s}_{(0)}(I,\mathbb{C})$ with
$$\| \Psi_{0} \|_{H^{s}_{(0)}} \leqslant \Gamma
\text{  and  }
|\langle \Psi_{0} , \phi_{1,\sigma} \rangle | \geqslant \gamma.$$
We claim that there exists $N=N(\Gamma,s,\epsilon,\gamma) \in \mathbb{N}^{*}$,
large enough, so that
\begin{equation} \label{cdt-Psi0-bis}
\begin{array}{lll}
\sum\limits_{k=N+1}^{\infty} |\langle \Psi_{0} , \phi_{k,\sigma} \rangle |^{2}
< \frac{\epsilon \gamma^{2}}{1-\epsilon}.
\end{array}
\end{equation}
Then, we consider the Lyapunov function
\begin{equation}\label{eq:lyap1}
\VVV (\Psi)=1-|\her{\Psi}{\phi_{1,\sigma}}|^2
-(1-\epsilon)\sum_{k=2}^N |\her{\Psi}{\phi_{k,\sigma}}|^2.
\end{equation}
Note that, this Lyapunov function depends on the constants
$\Gamma$, $s$, $\epsilon$, $\gamma$ through the choice of the cut-off dimension, $N$.
Just like~\cite{mirrahimi-cdc06}, it encodes two tasks: 1- it prevents
the $L^2$-mass lost through the high-energy eigenstates; 2- it
privileges the increase of the population in the first eigenstate.

When $\Psi$ solves $(\Sigma)$ with some control $u=\sigma + v$, we have
$$
\frac{d\VVV}{dt}= - 2 v(t) \Im \Big(
\sum_{k=1}^{N} a_{k}
\her{x\Psi}{\phi_{k,\sigma}}\her{\phi_{k,\sigma}}{\Psi}
\Big),
$$
where
\begin{equation} \label{def-a}
a_{1}:=1  \text{  and  } a_{k}:=1-\epsilon \text{  for  }
k=2,\cdots,N.
\end{equation}

Thus, the feedback law
\begin{equation}\label{eq:law1}
v(\Psi):=\varsigma \Im \Big( \sum_{k=1}^{N} a_{k}
\her{x\Psi}{\phi_{k,\sigma}}\her{\phi_{k,\sigma}}{\Psi} \Big),
\end{equation}
where $\varsigma>0$ is a positive constant, trivially ensures the
decrease of the Lyapunov function~\eqref{eq:lyap1}. We claim that,
the solution of (\ref{eq:main})-(\ref{eq:initial})-(\ref{eq:bords})
with initial condition $\Psi_{0}$ and control $u=\sigma + v(\Psi)$
satisfies
\begin{equation} \label{but-sigma}
\limsup\limits_{t \rightarrow + \infty}
\text{dist}_{L^{2}} (\Psi(t) , \mathcal{C}_{1,\sigma})^{2} \leqslant \epsilon.
\end{equation}
Note that, the claimed result here is much stronger than the one
provided in~\cite{mirrahimi-cdc06} for the finite potential well
problem. In fact, here, we claim the almost global approximate
stabilization of the system round the eigenstate $\phi_{1,\sigma}$.

The limit~\eqref{but-sigma} will be proved by studying the
$L^{2}(I,\CC)$-weak limits of $\Psi(t)$ when $t \rightarrow +
\infty$. Namely, let $(t_{n})_{n \in \mathbb{N}}$ be an increasing
sequence of positive real numbers such that $t_{n} \rightarrow +
\infty$. Since $\|\Psi(t_{n})\|_{L^{2}(I,\CC)} \equiv 1$, there
exists $\Psi_{\infty} \in L^{2}(I,\CC)$ such that, up to a
subsequence, $\Psi(t_{n}) \rightarrow \Psi_{\infty}$ weakly in
$L^{2}(I,\CC)$. Using the controllability of the linearized system
around $\phi_{1,\sigma}$ (which is equivalent to $\langle
\phi_{1,\sigma} , x \phi_{k,\sigma} \rangle \neq 0$ for every $k \in
\mathbb{N}^{*}$), we will be able to prove that $\Psi_{\infty}=\beta
\phi_{1,\sigma}$, where $\beta \in \CC$ and $|\beta|^{2} \geqslant
1-\epsilon$. This will imply~\eqref{but-sigma}.

Therefore, by weakening the stabilization property
(i.e. ask approximate stabilization instead of stabilization)
we avoid the compactness problem evoked at the begining of this section.

Note that, the controllability of the linearized system around the trajectory
$\phi_{1,\sigma}$ plays a crucial role here. This is why the
developed techniques  for $\sigma \neq 0$ can not be applied,
directly, to the case of $\sigma=0$.
\\

\begin{figure}[h]\psfrag{P}{$\Psi$}\psfrag{f}{$\phi_{1}$}\psfrag{s}{$\phi_{1,\sigma}$}
\begin{center}
  \includegraphics*[width=6cm]{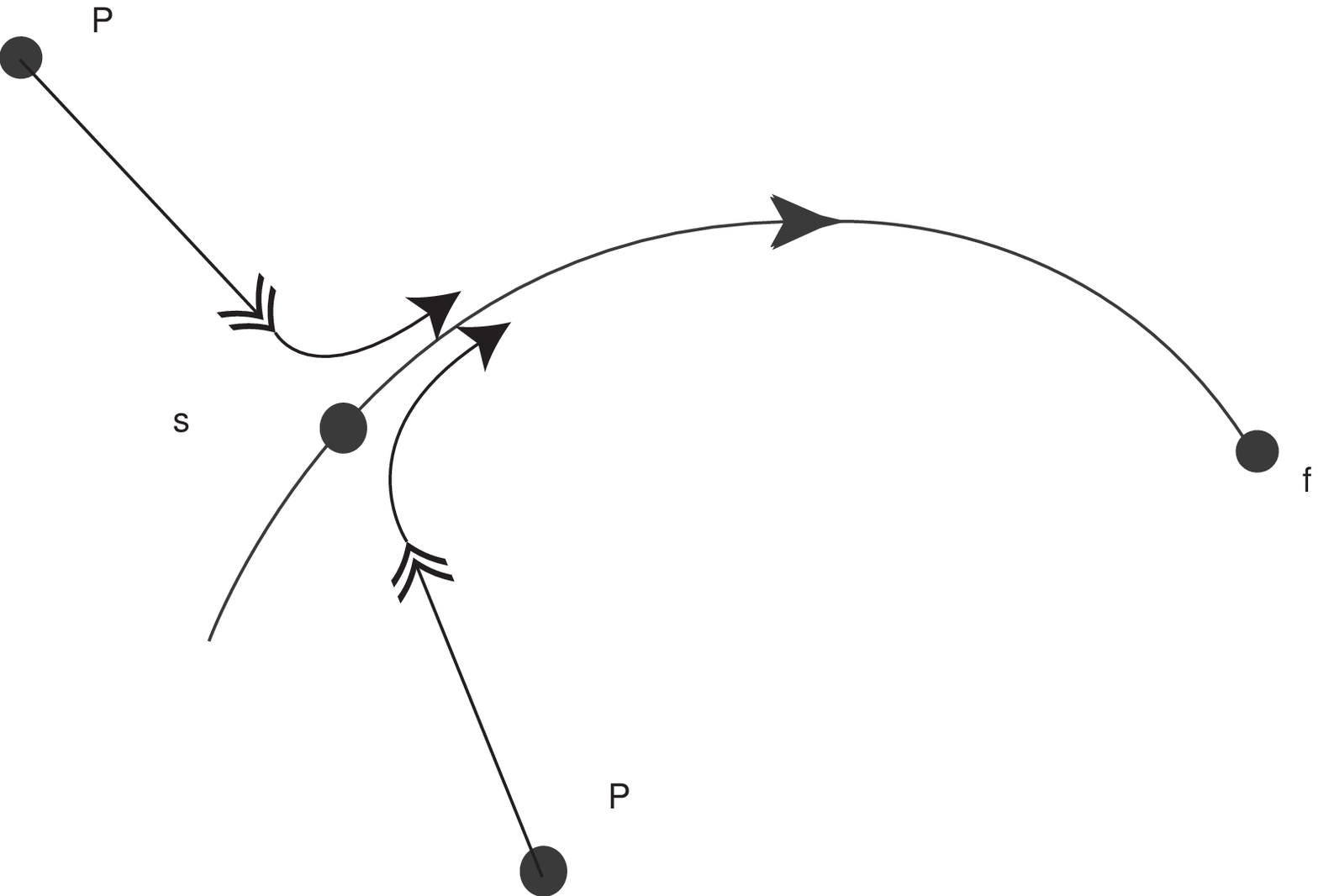}
  \caption{}
  \label{figaux}
  \end{center}
\end{figure}

Now, let us study the case $\sigma=0$. As emphasized above, the
previous strategy does not work for the approximate stabilization of
$\phi_{1}$ because the linearized system around $\phi_{1}$ is not
controllable. The idea is thus to use the above feedback
design~\eqref{eq:law1} with a dynamic $\sigma=\sigma(t)$ that
converges to zero as $t \rightarrow + \infty$. Formally, the
convergence of $\Psi$ toward $\CCC_{1,\sigma(t)}$ must happen at a
faster rate than that of $\sigma$ toward zero (see
Figure~\ref{figaux}).

In this aim, we consider the Lyapunov function
\begin{equation}\label{eq:lyap2}
\VVV(\Psi)=1-(1-\epsilon)\sum_{k=1}^N
|\her{\Psi}{\phi_{k,\sigma(\Psi)}}|^2-\epsilon|\her{\Psi}{\phi_{1,\sigma(\Psi)}}|^2,
\end{equation}
where the function $\Psi\mapsto \sigma(\Psi)$ is implicitly defined
as below
\begin{equation}\label{eq:implicit}
\sigma(\Psi)=\theta\left( \VVV(\Psi) \right),
\end{equation}
for a slowly varying real function $\theta$. We claim that such a
function $\sigma(\Psi)$ exists. When $\Psi$ solves $(\Sigma)$, we
have
$$\begin{array}{lll}
\frac{d\VVV}{dt}=
&
- 2 v(\Psi) \Im \Big(
\sum\limits_{k=1}^{N} a_{k}
\her{x\Psi}{\phi_{k,\sigma(\Psi)}}\her{\phi_{k,\sigma(\Psi)}}{\Psi}
\Big)
\\ &
- \frac{d\sigma(\Psi)}{dt} 2 \Re \Big( \sum_{k=1}^N  a_{k} \langle
\Psi , \phi_{k,\sigma(\Psi)} \rangle \langle
\frac{d\phi_{k,\sigma(\Psi)}}{d\sigma} , \Psi \rangle \Big)
\end{array}$$
where $\Re$ denotes the real part of a complex number, $(a_{k})_{1
\leqslant k \leqslant N}$ is defined by (\ref{def-a}) and the
notation $\frac{d\phi_{k,\sigma(\Psi)}}{d\sigma}$ means the
derivative of the map $\sigma \mapsto \phi_{k,\sigma}$ taken at the
point $\sigma=\sigma(\Psi)$. By definition of $\sigma(\Psi)$, we
have
$$\frac{d\sigma(\Psi)}{dt} = \theta'(\VVV(\Psi)) \frac{d\VVV}{dt}.$$
Thus, the feedback law $u(\Psi):=\sigma(\Psi)+v(\Psi)$ where
$$v(\Psi) := \varsigma \Im \Big(
\sum\limits_{k=1}^{N}  a_{k}
\her{x\Psi}{\phi_{k,\sigma(\Psi)}}\her{\phi_{k,\sigma(\Psi)}}{\Psi}
\Big)
$$
with $\varsigma>0$, ensures
$$\frac{d\VVV}{dt}= - 2 \varsigma \mu v(\Psi)^{2},$$
where
$$
\frac{1}{\mu}=1+2 \theta'(\VVV(\Psi)) \Re \Big( \sum_{k=1}^N a_{k}
\langle \Psi , \phi_{k,\sigma(\Psi)} \rangle \langle
\frac{d\phi_{k,\sigma(\Psi)}}{d\sigma} , \Psi \rangle \Big)
$$
is a positive constant, when $\|\theta'\|_{L^{\infty}}$ is small enough.
Thus $t \mapsto \VVV(\Psi(t))$ is not increasing.

We claim that, the solution of
(\ref{eq:main})-(\ref{eq:initial})-(\ref{eq:bords})
with initial condition $\Psi_{0}$ and control
$u=\sigma(\Psi) + v(\Psi)$
satisfies
\begin{equation}
\limsup\limits_{t \rightarrow + \infty}
\text{dist}_{L^{2}} (\Psi(t) , \mathcal{C}_{1})^{2} \leqslant \epsilon.
\end{equation}
Again, this will be proved by studying the $L^{2}(I,\CC)$-weak limits of $\Psi(t)$
when $t \rightarrow + \infty$.

\subsection{Structure of the article}
The rest of the paper is organized as follows.

The Section~\ref{sigma non nul} is dedicated to the proof of the
Theorem~\ref{Stab-Hs-sigma} when $\sigma \neq 0$. We derive this
theorem as a consequence of a stronger result stated in Theorem
\ref{Thm-Main result}.

This theorem and a straightforward corollary
(Corollary~\ref{Cor-sigma non nul}), leading to the
Theorem~\ref{Stab-Hs-sigma} in the case $\sigma\neq 0$, will be
stated in Subsection~\ref{Main result}. The
Subsection~\ref{Preliminaries} is dedicated to some preliminary
study needed for the proof of the Theorem~\ref{Thm-Main result} and
the Corollary~\ref{Cor-sigma non nul}. The proofs will be detailed
in Subsection~\ref{Proof of main theorem}.

The Section~\ref{sigma=0} is devoted to the proof of the Theorem
\ref{Stab-Hs-sigma} , in the case $\sigma =0$. Again, this theorem
will be derived as a consequence of a stronger result stated in
Theorem~\ref{Thm-Main result-0}.

In Subsection~\ref{Main result-sigma=0}, we state a Proposition
(Proposition~\ref{Existence-sigma}) ensuring the existence of the
implicit function $\sigma=\sigma(\Psi)$. Then, we state the Theorem
\ref{Thm-Main result-0} and a straightforward corollary (Corollary
\ref{Cor-Main result-sigma=0}), leading to the
Theorem~\ref{Stab-Hs-sigma} in the case $\sigma=0$. A preliminary
study, in preparation of the proof of the Theorem~\ref{Thm-Main
result-0} and the Corollary~\ref{Cor-Main result-sigma=0}, will be
performed in Subsection~\ref{Preliminaries-bis}. The proofs will be
detailed in Subsection~\ref{subsection:Cor-Main result-sigma=0}.

Finally, in Section~\ref{sec:numeric}, we provide some numerical
simulations to check out the performance of the control design on a
rather hard test case.

\section{Stabilization of $\CCC_{1,\sigma}$ with $\sigma \neq 0$}
\label{sigma non nul}
\subsection{Main result}\label{Main result}
The main result of Section~\ref{sigma non nul} is the following
theorem.
\begin{theorem}\label{Thm-Main result}
Let $N \in \NN^{*}$. There exists
$\sigma^{\sharp}=\sigma^{\sharp}(N)>0$ such that, for every $\sigma
\in (-\sigma^{\sharp},\sigma^{\sharp})-\{0\}$, $\gamma \in (0,1)$,
$\epsilon >0$, and $\Psi_{0} \in \SSS$ verifying
\begin{equation} \label{Hyp-Main result}
\begin{array}{lll}
\sum\limits_{k=N+1}^{\infty} |\langle \Psi_{0} , \phi_{k,\sigma}|^{2}
< \frac{\epsilon \gamma^{2}}{1-\epsilon}
& \text{   and   } &
|\langle \Psi_{0} , \phi_{1,\sigma} \rangle| \geqslant \gamma,
\end{array}
\end{equation}
the Cauchy problem
(\ref{eq:main})-(\ref{eq:initial})-(\ref{eq:bords}) with
$u(t)=\sigma+v_{\sigma,N,\epsilon}(\Psi(t))$,
\begin{equation} \label{Def-v}
v_{\sigma,N,\epsilon}(\Psi):=- \Im \left(
(1-\epsilon)\sum_{k=1}^{N} \langle x \Psi , \phi_{k,\sigma} \rangle
\overline{\langle \Psi , \phi_{k,\sigma} \rangle }
+ \epsilon \langle x \Psi , \phi_{1,\sigma} \rangle
\overline{\langle \Psi , \phi_{1,\sigma} \rangle} \right)
\end{equation}
has a unique weak solution $\Psi$. Moreover, this solution satisfies
\begin{equation} \label{Main result-ccl}
\liminf\limits_{t \rightarrow + \infty}
|\langle \Psi(t) , \phi_{1,\sigma} \rangle |^{2} \geqslant 1-\epsilon.
\end{equation}
\end{theorem}
The Theorem~\ref{Thm-Main result} provides an almost global
approximate stabilization. Indeed, any initial condition $\Psi_{0}
\in \SSS$ such that $\langle \Psi_{0} , \phi_{1,\sigma} \rangle \neq
0$ can be approximately moved to $\mathcal{C}_{1,\sigma}$. Notice
that the regularity assumption $\Psi_{0} \in
H^{s}_{(\sigma)}(I,\CC)$ stated in Theorem~\ref{Stab-Hs-sigma} is
not necessary for this purpose. Indeed, the feedback law depends on
the initial state through the choice of the cut-off dimension $N$.

The following corollary states that
the quantity $N$ appearing in the feedback law
may be uniform when $\Psi_{0}$ is in a given
bounded subset of $H^{s}_{(\sigma)}(I,\CC)$.

\begin{corollary} \label{Cor-sigma non nul}
Let $s>0$, $\epsilon>0$, $\Gamma>0$ and $\gamma \in (0,1)$. There
exists $\sigma^{**}=\sigma^{**}(\Gamma,s,\epsilon,\gamma)>0$
and $N=N(\Gamma,s,\epsilon,\gamma) \in \mathbb{N}^{*}$ such that,
for every $\sigma \in (-\sigma^{**},\sigma^{**})-\{0\}$, and
$\Psi_{0} \in H^{s}_{(\sigma)}(I,\CC) \cap \SSS$ verifying
\begin{equation} \label{Hyp-Cor-Main result}
\| \Psi_{0}  \|_{H^{s}_{(\sigma)}} \leqslant \Gamma
\text{  and  }
|\langle \Psi_{0} , \phi_{1,\sigma} \rangle | \geqslant \gamma,
\end{equation}
the Cauchy problem
(\ref{eq:main})-(\ref{eq:initial})-(\ref{eq:bords}) with
$u=\sigma+v_{\sigma,N,\epsilon}(\Psi)$, has a unique weak solution
$\Psi$. Moreover, this solution satisfies (\ref{Main result-ccl}).
\end{corollary}

\begin{remark}
The Theorem~\ref{Stab-Hs-sigma} in the case $\sigma \neq 0$ is a
direct consequence of the previous corollary. The feedback law
mentionned in the Theorem~\ref{Stab-Hs-sigma} is explicitly given in
the Corollary~\ref{Cor-sigma non nul}.

Notice that, in the particular case $\sigma \neq 0$, the Corollary
\ref{Cor-sigma non nul} is slightly more general than the Theorem
\ref{Stab-Hs-sigma}. In fact, the assumption \textquotedblleft
$\Psi_{0} \in H^{2} \cap H^{1}_{0}(I,\CC)$\textquotedblright is not
needed as we deal with weak solutions instead of strong ones.
Trivially, this solution will be a strong solution for $\Psi_{0} \in
H^{2} \cap H^{1}_{0}(I,\CC)$.

In the case $\sigma = 0$, this will no longer be the case : we will
need solutions in $C^{1}(\mathbb{R},L^{2})$ (for which the
assumption $\Psi_{0} \in H^{2} \cap H^{1}_{0}(I,\CC)$ is needed, see
the Proposition~\ref{Existence-classic}).
\end{remark}

\subsection{Preliminaries}\label{Preliminaries}
This section is devoted to the preliminary results, that will be
applied in the proof of the Theorem~\ref{Thm-Main result}.
\subsubsection{Eigenvalues and eigenvectors of $A_{\sigma}$}
\hspace{2cm}\vspace{.5cm}
\begin{proposition} \label{Vap-analytic}
For every $k \in \mathbb{N}^{*}$,
the eigenvalue $\sigma \mapsto \lambda_{k,\sigma} \in \RR$ and
the eigenstate $\sigma \mapsto \phi_{k,\sigma} \in (H^{2} \cap H^{1}_{0})(I,\CC)$
are analytic functions of $\sigma \in \mathbb{R}$ around
$\sigma=0$ and the expansion
$\lambda_{k,\sigma}=\lambda_{k} + \sigma^{2} \lambda_{k}^{(2)}
+o(\sigma^{2})$
holds with
\begin{equation} \label{lambda-k-2-explicite}
\lambda_{k}^{(2)}=\frac{1}{24k^{2}}-\frac{5}{8\pi^{2}k^{4}}.
\end{equation}
There exists $\sigma^{*}>0$, $C^{*}>0$ such that,
for every $\sigma_{0},\sigma_{1} \in (-\sigma^{*},\sigma^{*})-\{0\}$,
for every $k \in \mathbb{N}^{*}$,
\begin{alignat}{5}
\langle x \phi_{1,\sigma_{0}} , \phi_{k,\sigma_{0}} \rangle &\neq 0,
\label{ps-non-nul}\\
|\lambda_{k,\sigma_{0}} - \lambda_{k}| &\leqslant
\frac{C^{*}\sigma^2}{k},  \label{vap-borne-0}\\
\Big\| \frac{d\phi_{k,\sigma_{0}}}{d\sigma} \Big\|_{L^{2}}
&\leqslant \frac{C^{*}}{k}, \label{vep-k-1}\\
\Big\| \frac{d\phi_{k,\sigma_{0}}}{d\sigma}  \Big\|_{H^{1}_{0}}
&\leqslant C^{*}, \label{vep-k-1-H1}\\
\| \phi_{k,\sigma_{0}} - \phi_{k,\sigma_{1}} \|_{L^{2}} & \leqslant
\frac{C^{*}|\sigma_{0}-\sigma_{1}|}{k}.\label{diff-vep-k}
\end{alignat}
\end{proposition}
In the previous proposition, the notation
$\frac{d\phi_{k,\sigma_{0}}}{d\sigma}$ means the derivative of the
map $\sigma \mapsto \phi_{k\sigma}$ taken at the point
$\sigma=\sigma_{0}$. In the same way, we will use the notation
$\frac{d \lambda_{k,\sigma_{0}}}{d\sigma}$ for the derivative of the
map $\sigma \mapsto \lambda_{k,\sigma}$ at $\sigma=\sigma_{0}$.
\\

\textit{Proof of Proposition~\ref{Vap-analytic} : } We consider the
family of self-adjoint operators $A_{\sigma}=A-\sigma x$ in the
space $(H^2\cap H^1_0)(I,\CC)$. In this Banach space, the operator
$x$ (as a multiplication operator) is relatively bounded with
respect to $A$ with relative bound 0 (in the sense
of~\cite{kato-book-66}, page 190). Therefore $A_{\sigma}$ is a
self-adjoint holomorphic family of type (A)
(see~\cite{kato-book-66}, page 375). Thus the eigenvalues and the
eigenstates of $A_{\sigma}$ are holomorphic functions of $\sigma$.

Thanks to the Rayleigh-Schr\"{o}dinger perturbation theory,
we compute the first terms of the expansions
$$\begin{array}{ll}
\lambda_{k,\sigma}=\lambda_{k} + \sigma \lambda_{k}^{(1)} +
\sigma^{2} \lambda_{k}^{(2)} + \cdots, & \phi_{k,\sigma}=\phi_{k}+
\sigma \phi_{k}^{(1)} + \sigma^{2} \phi_{k}^{(2)}+\cdots.
\end{array}$$
Considering the first and second order terms of the equalities
$A_{\sigma}\phi_{k,\sigma}=\lambda_{k,\sigma} \phi_{k,\sigma}$,
$\|\phi_{k,\sigma}\|_{L^{2}}^{2}=1$, we get
\begin{equation} \label{expansion-1}
\begin{array}{lll}
- \frac{d^{2}}{dx^{2}} \phi_{k}^{(1)} - x \phi_{k}=
\lambda_{k} \phi_{k}^{(1)}+\lambda_{k}^{(1)}\phi_{k},
& &
\langle \phi_{k}^{(1)} , \phi_{k} \rangle = 0,
\end{array}
\end{equation}
\begin{equation} \label{expansion-2}
\begin{array}{lll}
- \frac{d^{2}}{dx^{2}} \phi_{k}^{(2)} - x \phi_{k}^{(1)}=
\lambda_{k}
\phi_{k}^{(2)}+\lambda_{k}^{(1)}\phi_{k}^{(1)}+\lambda_{k}^{(2)}\phi_{k},
\quad 2\langle \phi_{k}^{(2)} , \phi_{k} \rangle +
\|\phi_{k}^{(1)}\|_{L^{2}}^{2} = 0.
\end{array}
\end{equation}
Taking the Hermitian product of the first equality of
(\ref{expansion-1}) with $\phi_{k}$ and applying the parity
properties of $\phi_{k}$, we get $\lambda_{k}^{(1)}=0$. Considering
the Hermitian product of the first equality of (\ref{expansion-1})
with $\phi_{j}$, we get
\begin{equation} \label{phi_k^1 explicit}
\phi_{k}^{(1)}=\sum_{ j\in \NN^{*}, P(j) \neq P(k) }
\frac{\langle x \phi_{j} , \phi_{k} \rangle}{\lambda_{j}-\lambda_{k}}
\phi_{j},
\end{equation}
where the sum is taken over $j \in \mathbb{N}^{*}$ such that the
parity of $j$ is different from the parity of $k$. Taking the
Hermitian product of the first equality of (\ref{expansion-2}) with
$\phi_{k}$ we get $\lambda_{k}^{(2)}=-\langle x \phi_{k}^{(1)} ,
\phi_{k} \rangle$. Using (\ref{phi_k^1 explicit}) and the explicit
expression of $\langle x \phi_{k} , \phi_{j} \rangle$ computed
thanks to (\ref{eq:spectre}), we get
\begin{equation} \label{eq:sum}
\lambda_{k}^{(2)}=\frac{2^7}{\pi^{4}} \sum_{j\in \NN^{*}, P(j)\neq P(k)}
\frac{k^{2}j^{2}}{(k^{2}-j^{2})^{5}}.
\end{equation}
In order to simplify the above sum, we decompose the fraction
$$
F(X):=\frac{X^2}{(X-q)^5(X+q)^5}
$$
in the form\small
\begin{multline*}
F(X)=\frac{1}{2^5q^3}\left(\frac{1}{(X-q)^5}
-\frac{1}{(X+q)^5}\right)-\frac{1}{2^6q^4}\left(\frac{1}{(X-q)^4}+\frac{1}{(X+q)^4}\right)\\
-\frac{1}{2^7q^5}\left(\frac{1}{(X-q)^3}
-\frac{1}{(X+q)^3}\right)+\frac{5}{2^8q^6}\left(\frac{1}{(X-q)^2}+\frac{1}{(X+q)^2}\right)\\
-\frac{5}{2^8q^7}\left(\frac{1}{(X-q)}-\frac{1}{(X+q)}\right).
\end{multline*}\normalsize
Inserting this relation in the sum (\ref{eq:sum}) and simplifying,
we find
\begin{equation} \label{lambda-2-sum}
\begin{array}{ll}
\lambda_{k}^{(2)} =
&
\frac{1}{\pi^{4}} \left(
\frac{5}{2k^{5}} S_{k}^{1} - \frac{5}{2k^{4}} S_{k}^{2}
+\frac{1}{k^{3}} S_{k}^{3} + \frac{2}{k^{2}} S_{k}^{4}
-\frac{4}{k} S_{k}^{5} \right),
\end{array}
\end{equation}
where
$$S_{k}^{a}:=\sum_{j\in \NN^{*},P(j)\neq P(k)} \left(
\frac{1}{(j-k)^{a}} + \frac{(-1)^{a}}{(j+k)^{a}} \right) \text{ for
} a=1,\cdots,5.$$ We apply, now, the following well-known relations
for the Riemann $\zeta$-function:
$$
\zeta(2)=\sum_{j=1}^\infty \frac{1}{j^2}=\frac{\pi^2}{6}\qquad
\text{and}\qquad\zeta(4)=\sum_{j=1}^\infty
\frac{1}{j^4}=\frac{\pi^4}{90}.
$$
These relations imply
$$
\sum_{k=-\infty}^\infty\frac{1}{(2j+1)^2}
=\frac{\pi^2}{4}\qquad\text{and}\qquad\sum_{k=-\infty}^\infty\frac{1}{(2j+1)^4}
=\frac{\pi^4}{48},
$$
thus
$$\begin{array}{l}
S_{k}^{a}=
\left\lbrace\begin{array}{ll}
\frac{1}{k^{a}} & \text{ when } k \text{ is odd,}\\
0 & \text{ when } k \text{ is even}
\end{array}\right.
\text{ for } a=1,3,5,
\end{array}$$
$$\begin{array}{ll}
S_{k}^{2}=
\left\lbrace\begin{array}{ll}
\frac{\pi^{2}}{4}-\frac{1}{k^{2}}
& \text{ when } k \text{ is odd,}\\
\frac{\pi^{2}}{4}
& \text{ when } k \text{ is even,}
\end{array}\right.
&
S_{k}^{4}=
\left\lbrace\begin{array}{ll}
\frac{\pi^{4}}{48}-\frac{1}{k^{4}}
& \text{ when } k \text{ is odd,}\\
\frac{\pi^{4}}{48}
& \text{ when } k \text{ is even.}
\end{array}\right.
\end{array}$$
Inserting this in (\ref{lambda-2-sum}), we get (\ref{lambda-k-2-explicite}).

The relation (\ref{ps-non-nul}) is proved in \cite[Proposition
1]{beauchard-04}. The bound (\ref{vap-borne-0}) is given in
\cite[Chapter 17 Example 2.14, Chapter 2 Problem 3.7]{kato-book-66}.
The inequality (\ref{vep-k-1}) is proved in \cite[Proposition
42]{beauchard-04}. The bound (\ref{vep-k-1-H1}) is a consequence of
(\ref{vep-k-1}). Indeed, considering the Hermitian product in
$L^{2}(I,\CC)$ of $\frac{d\phi_{k,\sigma_{0}}}{d\sigma}$ with the
equation
$$A_{\sigma_{0}} \frac{d\phi_{k,\sigma_{0}}}{d\sigma}
- x \phi_{k,\sigma_{0}} =
\lambda_{k,\sigma_{0}} \frac{d\phi_{k,\sigma_{0}}}{d\sigma}
+ \frac{d\lambda_{k,\sigma_{0}}}{d\sigma}  \phi_{k,\sigma_{0}},$$
and using (\ref{vep-k-1}) together with
the orthogonality between $\phi_{k,\sigma_{0}}$
and $\frac{d\phi_{k,\sigma_{0}}}{d\sigma}$
(which is a consequence of $\|\phi_{k,\sigma}\|_{L^{2}}^{2} \equiv 1$),
we get
$$-\frac{1}{2}
\Big\|  \frac{d\phi_{k,\sigma_{0}}}{d\sigma}  \|_{H^{1}_{0}}^{2}
\leqslant |\sigma_{0}| \left( \frac{C^{*}}{k} \right)^{2} +
\frac{C^{*}}{k} + \left( \frac{\pi^{2}k^{2}}{2} + C^{*}
\sigma_{0}^{2} \right) \left( \frac{C^{*}}{k} \right)^{2},$$ which
gives (\ref{vep-k-1-H1}). Finally, (\ref{diff-vep-k}) is a
consequence of (\ref{vep-k-1}). \endproof

\begin{proposition} \label{VAP}
Let $N \in \NN^{*}$. There exists
$\sigma^{\sharp}=\sigma^{\sharp}(N)>0$ such that, for every $\sigma
\in (-\sigma^{\sharp},\sigma^{\sharp})-\{0\}$, $j_{2}, k_{2} \in
\NN^{*}$, and $j_{1},k_{1} \in \{1,\cdots,N\}$, verifying $j_{1}
\neq j_{2}$ and $k_{1} \neq k_{2}$,
\begin{equation} \label{Hyp-k1-k2-j1-j2}
\lambda_{k_{1},\sigma} - \lambda_{k_{2},\sigma}
=\lambda_{j_{1},\sigma} - \lambda_{j_{2},\sigma}
\end{equation}
implies $(j_{1},j_{2})=(k_{1},k_{2})$.
\end{proposition}

\textit{Proof of Proposition~\ref{VAP}: } Let $C^{*}$ be as in
Proposition~\ref{Vap-analytic} and $\sigma \in
(-\sigma^{\sharp}_{0},\sigma^{\sharp}_{0})$ where
\begin{equation} \label{Def-sigma-sharp}
\sigma^{\sharp}_{0}:= \frac{\pi}{4\sqrt{C^{*}}}.
\end{equation}
First, we prove the equality (\ref{Hyp-k1-k2-j1-j2}) to be
impossible when $j_{2} \neq k_{2}$ and
\begin{equation} \label{hyp-absurd-1}
\max\{ j_{2} , k_{2} \} > \frac{N^{2}+1}{2}.
\end{equation}
We argue by contradiction. Let us assume the existence of $j_{2},
k_{2} \in \NN^{*}$, $j_{1},k_{1} \in \{1,\cdots,N\}$, with $j_{1}
\neq j_{2}$, $k_{1} \neq k_{2}$, $j_{2} \neq k_{2}$, such that
(\ref{hyp-absurd-1}) and (\ref{Hyp-k1-k2-j1-j2}) hold. Without loss
of generality, we may assume that $\max\{ j_{2} , k_{2} \}=j_{2} >
\frac{N^{2}+1}{2}$. Using (\ref{vap-borne-0}), we get
$$\begin{array}{ll}
\lambda_{j_{2},\sigma}-\lambda_{k_{2},\sigma}
& \geqslant \frac{\pi^{2}}{2} (j_{2}^{2} - k_{2}^{2}) - 2 C^{*} \sigma^{2}
\\& \geqslant \frac{\pi^{2}}{2} (j_{2}^{2} - (j_{2}-1)^{2}) - 2 C^{*} \sigma^{2}
\\& \geqslant \frac{\pi^{2}}{2} (2j_{2}-1) - 2 C^{*} \sigma^{2},
\end{array}$$
$$\begin{array}{ll}
\lambda_{j_{1},\sigma} - \lambda_{k_{1},\sigma}
& \leqslant \frac{\pi^{2}}{2} (N^{2}-1) + 2C^{*}\sigma^{2}.
\end{array}$$
Using the equality of the left hand sides of these
inequalities, together with (\ref{Def-sigma-sharp}), we get
$$j_{2} \leqslant \frac{N^{2}}{2} + \frac{8C^{*}\sigma_{0}^{\sharp 2}}{\pi^{2}}
\leqslant \frac{N^{2}+1}{2} ,$$
which is a contradiction.

Therefore, it is sufficient to prove the Proposition~\ref{VAP} for
$j_2,k_2\in\{1,\cdots,[(N^2+1)/2]\}$. Moreover, it is sufficient to
prove that, for every $j_1,k_1\in\{1,\cdots,N\}$ and
$j_2,k_2\in\{1,\cdots,[(N^2+1)/2]\}$, with $j_{1} \neq j_{2}$,
$k_{1} \neq k_{2}$, $(j_{1},j_{2}) \neq (k_{1},k_{2})$, there exists
$\sigma^{\sharp}_{j_{1},k_{1},j_{2},k_{2}} \in
(0,\sigma^{\sharp}_{0})$ such that, for every $\sigma \in
(-\sigma^{\sharp}_{j_{1},k_{1},j_{2},k_{2}},
\sigma^{\sharp}_{j_{1},k_{1},j_{2},k_{2}})$, (\ref{Hyp-k1-k2-j1-j2})
does not hold. Indeed, then, the following choice of
$\sigma^\sharp(N)$ concludes the proof of the Proposition~\ref{VAP},
$$\begin{array}{ll}
\sigma^{\sharp}(N):=\min \{
\sigma^{\sharp}_{j_{1},k_{1},j_{2},k_{2}}; &
j_1,k_1\in\{1,\cdots,N\}, j_2,k_2\in\{1,\cdots,(N^2+1)/2\},
\\ &
(j_{1},j_{2}) \neq (k_{1},k_{2}),
j_{1} \neq j_{2},
k_{1} \neq k_{2} \}.
\end{array}$$
Let $j_1,k_1\in\{1,\cdots,N\}$, $j_2,k_2\in\{1,\cdots,(N^2+1)/2\}$
be such that $j_{1} \neq j_{2}$, $k_{1} \neq k_{2}$, $(j_{1},j_{2})
\neq (k_{1},k_{2})$. We argue by contradiction. Let us assume that,
for every $\sigma^{\sharp}_{1}>0$, there exists $\sigma \in
(-\sigma^{\sharp}_{1},\sigma^{\sharp}_{1})$ such that
(\ref{Hyp-k1-k2-j1-j2}) holds. Using the analyticity of both sides
in (\ref{Hyp-k1-k2-j1-j2}) with respect to $\sigma$, at $\sigma=0$,
this assumption implies that
$$\lambda_{k_{1}}^{(2)} - \lambda_{k_{2}}^{(2)}
=\lambda_{j_{1}}^{(2)} - \lambda_{j_{2}}^{(2)}.$$
Using (\ref{lambda-k-2-explicite}) together with a rationality argument, we get
$$\begin{array}{ll}
\frac{1}{k_{1}^{2}} - \frac{1}{k_{2}^{2}} = \frac{1}{j_{1}^{2}} - \frac{1}{j_{2}^{2}},
&
\frac{1}{k_{1}^{4}} - \frac{1}{k_{2}^{4}} = \frac{1}{j_{1}^{4}} - \frac{1}{j_{2}^{4}}.
\end{array}$$
Since $k_{1} \neq k_{2}$ and $j_{1} \neq j_{2}$,
we deduce from the previous equalities that
$$\begin{array}{ll}
\frac{1}{k_{1}^{2}} - \frac{1}{k_{2}^{2}} = \frac{1}{j_{1}^{2}} - \frac{1}{j_{2}^{2}},
&
\frac{1}{k_{1}^{2}} + \frac{1}{k_{2}^{2}} = \frac{1}{j_{1}^{2}} + \frac{1}{j_{2}^{2}}.
\end{array}$$
Therefore $k_{1}=j_{1}$ and $k_{2}=j_{2}$, which is a contradiction.
\endproof

\subsubsection{Solutions of the Cauchy Problem}
\hspace{2cm}\vspace{.5cm}
\begin{proposition} \label{Existence}
Let $\sigma \in \mathbb{R}$, $N \in \NN^{*}$, $\epsilon>0$. For
every $\Psi_{0} \in \SSS$, there exists a unique weak solution
$\Psi$ of (\ref{eq:main})-(\ref{eq:initial})-(\ref{eq:bords}) with
$u(t)=\sigma+v_{\sigma,N,\epsilon} (\Psi(t))$, i.e. $\Psi \in
C^{0}(\RR,\SSS) \cap C^{1}(\RR,H^{-2}_{(0)}(I,\CC))$, the equality
(\ref{eq:main}) holds in $H^{-2}_{(0)}(I,\CC)$ for every $t \in \RR$
and the equality (\ref{eq:initial}) holds in $\SSS$.
\end{proposition}

\textit{Proof of Proposition~\ref{Existence} : } Let $\sigma \in
\mathbb{R}$, $N \in \NN$, $\epsilon>0$, $\Psi_{0} \in \SSS$ and
$T>0$ be such that
\begin{equation} \label{Hyp-contraction}
TNe^{NT} < 1.
\end{equation}
In order to build solutions on $[0,T]$,
we apply the Banach fixed point theorem to the map
$$\begin{array}{cccc}
\Theta : & C^{0}([0,T],\SSS)    & \rightarrow & C^{0}([0,T],\SSS) \\
         & \xi                  & \mapsto     & \Psi
\end{array}$$
where $\Psi$ is the solution of
(\ref{eq:main})-(\ref{eq:initial})-(\ref{eq:bords}) with
$u(t)=\sigma + v_{\sigma,N,\epsilon}(\xi(t))$.

The map $\Theta$ is well defined and maps $C^{0}([0,T],\SSS)$ into
itself. Indeed, when $\xi \in C^{0}([0,T],\SSS)$, $u:t \mapsto
\sigma + v_{\sigma,N,\epsilon}(\xi(t))$ is continuous and thus the
Proposition~\ref{Existence-classic} ensures the existence of a
unique weak solution $\Psi$. Notice that the map $\Theta$ takes
values in $C^{0}([0,T],\SSS)) \cap C^{1}([0,T],H^{-2}_{(0)})$.

Let us prove that $\Theta$ is a contraction of $C^{0}([0,T],\SSS)$.
Let $\xi_{j} \in C^{0}([0,T],\SSS)$,
$v_{j}:=v_{\sigma,N,\epsilon}(\xi_{j})$,
$\Psi_{j}:=\Theta(\xi_{j})$, for $j=1,2$
and $\Delta:=\Psi_{1}-\Psi_{2}$. We have
$$\Delta(t)=i \int_{0}^{t}
e^{-iA_{\sigma}(t-s)} [ v_{1} x \Delta (s) + (v_{1}-v_{2}) x
\Psi_{2}(s) ]ds.$$ Thanks to (\ref{Def-v}), we have
$\|v_{j}\|_{L^{\infty}(0,T)} \leqslant N$ for $j=1,2$ and
$\|v_{1}-v_{2}\|_{L^{\infty}(0,T)} \leqslant 2N
\|\xi_{1}-\xi_{2}\|_{C^{0}([0,T],L^{2})}$. Thus
\begin{equation} \label{Delta-intermediaire}
\|\Delta(t)\|_{L^{2}} \leqslant
\int_{0}^{t}
N \|\Delta(s)\|_{L^{2}}
+ N \|\xi_{1}-\xi_{2}\|_{C^{0}([0,T],L^{2})} ds.
\end{equation}
Therefore, the Gronwall Lemma implies
$$\|\Delta(t)\|_{C^{0}([0,T],L^{2})} \leqslant
\|\xi_{1}-\xi_{2}\|_{C^{0}([0,T],L^{2})} NT e^{NT},$$ and so
(\ref{Hyp-contraction}) ensures that $\Theta$ is a contraction of
the Banach space $C^{0}([0,T],\SSS)$. Therefore, there exists a
fixed point $\Psi \in C^{0}([0,T],\SSS)$ such that
$\Theta(\Psi)=\Psi$. Since $\Theta$ takes values in
$C^{0}([0,T],\SSS) \cap C^{1}([0,T],H^{-2}_{(0)})$, necessarily
$\Psi$ belongs to this space, thus, it is a weak solution of
(\ref{eq:main})-(\ref{eq:initial})-(\ref{eq:bords}) on $[0,T]$.

Finally, we have introduced a time $T>0$ and, for every $\Psi_{0}
\in \SSS$, we have built a weak solution $\Psi \in
C^{0}([0,T],\SSS)$ of
(\ref{eq:main})-(\ref{eq:initial})-(\ref{eq:bords}) on $[0,T]$.
Thus, for a given initial condition $\Psi_{0} \in \SSS$, we can
apply this result on $[0,T]$, $[T,2T]$, $[2T,3T]$ etc. This proves
the existence and uniqueness of a global solution for the
closed-loop system. \endproof

\begin{proposition} \label{Convergence-H^{-1}}
Let $\sigma>0$, $N \in \NN$, $\epsilon >0$,
$(\Psi_{0}^{n})_{n \in \NN}$ be a sequence of $\SSS$
and $\Psi_{0}^{\infty} \in L^{2}$ with
$\|\Psi_{0}^{\infty}\|_{L^{2}} \leqslant 1$
be such that
$$\lim\limits_{n \rightarrow + \infty} \Psi_{0}^{n} = \Psi_{0}^{\infty}
\text{ strongly in } H^{-1}(I,\CC).$$ Let $\Psi^{n}$ (resp.
$\Psi^{\infty}$) be the weak solution of
(\ref{eq:main})-(\ref{eq:initial})-(\ref{eq:bords}) with
$u(t)=\sigma+v_{\sigma,N,\epsilon}(\Psi^{n})$ (resp. with
$u(t)=\sigma+v_{\sigma,N,\epsilon}(\Psi^{\infty}(t))$). Then, for
every $\tau>0$,
$$\lim\limits_{n \rightarrow + \infty} \Psi^{n}(\tau) = \Psi^{\infty}(\tau)
\text{ strongly in } H^{-1}(I,\CC).$$
\end{proposition}

\textit{Proof of Proposition~\ref{Convergence-H^{-1}} : } Let us
recall that the space $H^{-1}(I,\CC)$ (dual space of
$H^{1}_{0}(I,\CC)$ for the $L^{2}(I,\CC)$-Hermitian product)
coincides with $H^{-1}_{(0)}(I,\CC)$ and that
$\sqrt{2}\|.\|_{H^{-1}}=\|.\|_{H^{-1}_{(0)}}$ (because
$\|.\|_{H^{1}_{0}}= \sqrt{2} \|.\|_{H^{1}_{(0)}}$). We introduce
$\mathcal{C}>0$ such that,
\begin{equation} \label{Def-cst-x}
\| x \varphi \|_{H^{-1}} \leqslant \mathcal{C} \| \varphi
\|_{H^{-1}} \text{  ,  } \forall \varphi \in H^{-1}(I,\CC).
\end{equation}
Such a constant does exist. Indeed, for every $\xi \in
H^{1}_{0}(I,\CC)$, $x\xi \in H^{1}_{0}(I,\CC)$ and
$$\| x\xi \|_{H^{1}_{0}}=
\left( \int_{I} | x \xi' + \xi |^{2} dx \right)^{1/2}
\leqslant
\|\xi'\|_{L^{2}} (1+C_{P})$$
where $C_{P}$ is the Poincar\'e constant on $I$. Thus, for
$\varphi \in H^{-1}(I,\CC)$, we have
$$\begin{array}{ll}
\|x\varphi\|_{H^{-1}(I,\CC)}
& =
\sup \left\{  \langle x \varphi , \xi \rangle ;
\xi \in H^{1}_{0}(I,\CC), \|\xi\|_{H^{1}_{0}}=1 \right\}
\\ & \leqslant
\sup \left\{ \|\varphi\|_{H^{-1}} \|x\xi\|_{H^{1}_{0}};
\xi \in H^{1}_{0}(I,\CC), \|\xi\|_{H^{1}_{0}}=1 \right\}
\\ & \leqslant (1+C_{P})\|\xi\|_{H^{-1}}.
\end{array}$$

In order to simplify the notations, in this proof, we write
$v(\Psi)$ instead of $v_{\sigma,N,\epsilon}(\Psi)$.
We have
$$\begin{array}{lll}
(\Psi^{n}-\Psi^{\infty})(t)
& = &
e^{-iA t}(\Psi^{n}_{0}-\Psi^{\infty}_{0})
+i \int_{0}^{t} e^{-iA(t-s)} \sigma x (\Psi^{n}-\Psi^{\infty})(s) ds
\\ & &
+i \int_{0}^{t} e^{-iA(t-s)}
[v(\Psi^{n}(s))-v(\Psi^{\infty}(s))]x\Psi^{n}(s) ds
\\ & &
+i \int_{0}^{t} e^{-iA(t-s)}
v(\Psi^{\infty}(s)) x [\Psi^{n}(s)-\Psi^{\infty}(s)] ds.
\end{array}$$
Using (\ref{Def-v}), $\|\Psi^{n}(s)\|_{L^{2}}=1$,
$\|\Psi^{\infty}(s)\|_{L^{2}} \leqslant 1$ and the fact that
$\phi_{k,\sigma}, x \phi_{k,\sigma} \in H^{1}_{0}(I,\CC)$ for
$k=1,\cdots,N$, we get
\begin{equation} \label{v-difference}
|v(\Psi^{n}(s))-v(\Psi^{\infty}(s))| \leqslant 2N
\mathcal{C} C_{\sigma}(N) \| (\Psi^{n}-\Psi^{\infty})(s) \|_{H^{-1}},
\end{equation}
where $C_{\sigma}(N):=\sup\{ \|
\phi_{k,\sigma}\|_{H^{1}_{0}(I,\CC)}; k \in \{1,\cdots,N\} \}$. The
semigroup $e^{-iAt}$ preserves the $H^{-1}$-norm thus, using $|
v(\Psi^{\infty}(s)) | \leqslant N$ and (\ref{v-difference}), we get
$$\begin{array}{lll}
\|(\Psi^{n}-\Psi^{\infty})(t)\|_{H^{-1}} & & \leqslant
\|\Psi^{n}_{0}-\Psi^{\infty}_{0}\|_{H^{-1}}
\\ & &
+ \mathcal{C} \int_{0}^{t} ( |\sigma| + 2NC_{\sigma}(N) + N )
\|\Psi^{n}(s)-\Psi^{\infty}(s)\|_{H^{-1}} ds.
\end{array}$$
We conclude thanks to the Gronwall Lemma. \endproof

\subsection{Proof of Theorem~\ref{Thm-Main result} and Corollary~\ref{Cor-sigma non nul}}
\label{Proof of main theorem} \hspace{2 cm}\vspace{.5 cm}

\textit{Proof of Theorem~\ref{Thm-Main result} :} Let $N \in
\NN^{*}$. Let $\sigma^{*}>0$ be as in Proposition~\ref{Vap-analytic}
and $\sigma^{\sharp}=\sigma^{\sharp}(N)$ be as in Proposition
\ref{VAP}. Let $\sigma^{**}:=\min\{ \sigma^{*} , \sigma^{\sharp}
\}$.

Let $\sigma \in (-\sigma^{**},\sigma^{**})-\{0\}$, $\gamma \in
(0,1)$, $\epsilon >0$, $\Psi_{0} \in \SSS$ with (\ref{Hyp-Main
result}) and $\Psi$ be the weak solution of
(\ref{eq:main})-(\ref{eq:initial})-(\ref{eq:bords}) with
$u(t)=\sigma+v_{\sigma,N,\epsilon}(\Psi(t))$ given by
Proposition~\ref{Existence}. For $\varphi \in L^{2}(I,\CC)$, we
define
\begin{equation} \label{Def-V-Lyap}
\mathcal{V}_{\sigma, N, \epsilon}(\varphi):= 1-|\langle \varphi ,
\phi_{1,\sigma}\rangle |^{2} - (1-\epsilon) \sum_{k=2}^{N} |\langle
\varphi , \phi_{k,\sigma} \rangle |^{2}.
\end{equation}
Since $\Psi \in C^{1}(\mathbb{R},H^{-2}_{(0)}(I,\CC))$ and
$\phi_{k,\sigma} \in H^{2}_{(0)}(I,\CC)$, $t \mapsto \mathcal{V}_{N,
\sigma, \epsilon}(\Psi(t))$ is $C^{1}$. Using (\ref{eq:main}),
integrations by parts and $a_{1}:=1$, $a_{k}:=1-\epsilon$ when $k
\geqslant 2$, we get
\begin{equation} \label{Derivee-Lyapounov}
\begin{array}{ll}
\frac{d}{dt} \mathcal{V}_{\sigma,N,\epsilon}(\Psi)
&
= - 2  \Re
\left( \sum_{k=1}^{N} a_{k} \langle - i A_\sigma \Psi + i
v_{\sigma,N,\epsilon}(\Psi) x \Psi , \phi_{k,\sigma} \rangle
\overline{ \langle \Psi , \phi_{k,\sigma} \rangle} \right),
\\ &
=-2v_{\sigma,N,\epsilon}(\Psi(t))^{2}.
\end{array}
\end{equation}
Thus, $t \mapsto \mathcal{V}_{\sigma, N, \epsilon}(\Psi(t))$ is a
non increasing function. There exists $\alpha \in
[0,\mathcal{V}_{\sigma, N, \epsilon}(\Psi_{0})]$ such that
$\mathcal{V}_{\sigma, N, \epsilon}(\Psi(t)) \rightarrow \alpha$ when
$t \rightarrow + \infty$. Since $\Psi_{0} \in \SSS$ and
(\ref{Hyp-Main result}) holds we have
$$\begin{array}{ll}
\mathcal{V}_{\sigma, N, \epsilon}(\Psi_{0}) & = 1-(1-\epsilon)
\sum_{k=1}^{N} |\langle \Psi , \phi_{k,\sigma} \rangle |^{2} -
\epsilon |\langle \Psi , \phi_{1,\sigma}\rangle |^{2}
\\ & =
1-(1-\epsilon)\left( 1 -
\sum_{k=N+1}^{\infty} |\langle \Psi , \phi_{k,\sigma} \rangle |^{2} \right)
- \epsilon |\langle \Psi , \phi_{1,\sigma}\rangle |^{2}
\\ & <
1 - (1-\epsilon) \left( 1 - \frac{\epsilon \gamma^{2}}{1-\epsilon} \right)
- \epsilon \gamma^{2}
\\ & < \epsilon,
\end{array}$$
thus $\alpha \in [0,\epsilon)$.

Let $(t_{n})_{n \in \mathbb{N}}$ be an increasing sequence
of positive real numbers such that $t_{n} \rightarrow + \infty$
when $n \rightarrow + \infty$. Since $\|\Psi(t_{n})\|_{L^{2}}=1$
for every $n \in \mathbb{N}$, there exists
$\Psi_{\infty} \in L^{2}(I,\CC)$ such that, up to an extraction
$$\Psi(t_{n}) \rightarrow \Psi_{\infty}
\text{ weakly in } L^{2}(I,\CC)
\text{ and strongly in } H^{-1}(I,\CC).$$
Let $\xi$ be the solution of
$$\left\lbrace \begin{array}{l}
i \frac{\partial \xi}{\partial t} = A_{\sigma} \xi -
v_{\sigma,N,\epsilon}(\xi(t)) x \xi
\text{ , } x \in I \text{  ,  } t \in (0,+\infty),\\
\xi(t,\pm 1/2)=0,\\
\xi(0)=\Psi_{\infty}.
\end{array} \right.$$
Thanks to the Proposition~\ref{Convergence-H^{-1}}, for every
$\tau>0$, $\Psi(t_{n}+\tau) \rightarrow \xi(\tau)$ strongly in
$H^{-1}(I,\mathbb{C})$ when $n \rightarrow + \infty$. Thus
$\mathcal{V}_{\sigma, N, \epsilon}(\Psi(t_{n}+\tau)) \rightarrow
\mathcal{V}_{\sigma, N, \epsilon}(\xi(\tau))$ when $n \rightarrow +
\infty$, because $\mathcal{V}_{\sigma,N,\epsilon}(.)$ is continuous
for the $L^2$-weak topology. Therefore $\mathcal{V}_{\sigma, N,
\epsilon}(\xi(\tau)) \equiv \alpha$. Furthermore, the relation
(\ref{Derivee-Lyapounov}) holds when $\Psi$ is replaced by $\xi$, and
thus $v_{\sigma,N,\epsilon}(\xi(\tau)) \equiv 0$ and $\xi$ solves
$$\left\lbrace \begin{array}{l}
i \frac{\partial \xi}{\partial t} = A_{\sigma} \xi
\text{ , } x \in I \text{  ,  } t \in (0,+\infty),\\
\xi(t,\pm 1/2)=0,\\
\xi(0)=\Psi_{\infty}.
\end{array} \right.$$
Therefore, we have
$$\xi(\tau)=\sum_{k=1}^{\infty} \langle \Psi_{\infty} , \phi_{k,\sigma} \rangle
\phi_{k,\sigma} e^{-i\lambda_{k,\sigma}\tau}.$$ The equality
$v_{\sigma,N,\epsilon}(\xi) \equiv 0$, then, gives
\begin{equation} \label{feedback=0}
\begin{array}{ll}
\Im  \left(
\sum\limits_{k=1}^{N} \sum\limits_{j \in \NN^{*},j \neq k}
a_{k}
\langle \Psi_{\infty} , \phi_{j,\sigma} \rangle
\langle x \phi_{j,\sigma} , \phi_{k,\sigma} \rangle
\overline{ \langle \Psi_{\infty} , \phi_{k,\sigma} \rangle }
e^{i(\lambda_{k,\sigma}-\lambda_{j,\sigma})\tau}
\right) \equiv 0.
\end{array}
\end{equation}
Let
$\omega_{(k_{1},k_{2})}:=\lambda_{k_{1},\sigma}-\lambda_{k_{2},\sigma}$
for every $k_{1}, k_{2} \in \NN^{*}$ and $\mathcal{S}:=\{
(k_{1},k_{2}) ; k_{1} \in \{1,\cdots,N\},$ $k_{2} \in \NN^{*}, k_{1}
\neq k_{2} \}$. Thanks to the Proposition~\ref{VAP}, all the
frequencies $\omega_{K}$ for $K \in \mathcal{S}$ are different.
Moreover, there exists a uniform gap $\delta >0$ such that, for
every $\omega, \tilde{\omega} \in \{ \pm \omega_{K} ; K \in
\mathcal{S} \}$ with $\omega \neq \tilde{\omega}$, then $|\omega-
\tilde{\omega}| \geqslant \delta$. Thus, for $T>0$ large enough,
there exists $C=C(T)>0$ such that the Ingham inequality
$$ \sum\limits_{ K \in \mathcal{S}} |a_{K}|^{2}
\leqslant C \int_{0}^{T} \Big| \sum\limits_{K \in \mathcal{S}} a_{K}
e^{i\omega_{K} t} \Big|^{2} dt$$ holds, for every $(a_{K})_{K \in
\mathcal{S}} \in l^{2}(\mathcal{S},\CC)$ (see \cite[Theorem
1.2.9]{Krabs}). The equality (\ref{feedback=0}) implies, in
particular,
$$\langle \Psi_{\infty} , \phi_{j,\sigma} \rangle
\langle x \phi_{j,\sigma} , \phi_{1,\sigma} \rangle
\overline{ \langle \Psi_{\infty} , \phi_{1,\sigma} \rangle }=0,
\forall j \geqslant 2.$$
Thanks to (\ref{ps-non-nul}), we get
\begin{equation} \label{coeff=0}
\langle \Psi_{\infty} , \phi_{j,\sigma} \rangle
\overline{ \langle \Psi_{\infty} , \phi_{1,\sigma} \rangle }=0,
\forall j \geqslant 2.
\end{equation}
Let us prove that
\begin{equation} \label{Mode-1}
\langle \Psi_{\infty} , \phi_{1,\sigma} \rangle \neq 0.
\end{equation}
Since $\|\Psi^{\infty}\|_{L^{2}} \leqslant 1$, we have
$$\begin{array}{ll}
\mathcal{V}_{\sigma, N, \epsilon}(\Psi_{\infty}) & \geqslant 1-|
\langle \Psi^{\infty} , \phi_{1,\sigma} \rangle |^{2} - (1-\epsilon)
\sum_{k=2}^{\infty} | \langle \Psi^{\infty} , \phi_{k,\sigma}
\rangle |^{2}
\\ &
=
1-| \langle \Psi^{\infty} , \phi_{1,\sigma} \rangle |^{2}
- (1-\epsilon)[ \|\Psi^{\infty}\|_{L^{2}}^{2} -
| \langle \Psi^{\infty} , \phi_{1,\sigma} \rangle |^{2} ]
\\ &
\geqslant
\epsilon - \epsilon | \langle \Psi^{\infty} , \phi_{1,\sigma} \rangle |^{2}.
\end{array}$$
Moreover, $\mathcal{V}_{\sigma, N,
\epsilon}(\Psi_{\infty})=\alpha<\epsilon$, thus
$$\epsilon > \epsilon - \epsilon|\langle \Psi_{\infty} , \phi_{1,\sigma}|^{2},$$
which gives (\ref{Mode-1}). Therefore (\ref{coeff=0}) justifies the
existence of $\beta \in \CC$ with $|\beta| \leqslant 1$ such that
$\Psi_{\infty}=\beta \phi_{1,\sigma}$. Then, $\epsilon >
\alpha=\mathcal{V}_{N, \sigma, \epsilon}(\Psi_{\infty})
=1-|\beta|^{2}$, thus $|\beta|^{2}>1-\epsilon$. Finally, we have
$$\lim\limits_{n \rightarrow + \infty}
|\langle \Psi(t_{n}),\phi_{1,\sigma} \rangle|^{2} = |\langle
\Psi_{\infty},\phi_{1,\sigma} \rangle|^{2} =
|\beta|^{2}>1-\epsilon.$$ This holds for every sequence $(t_{n})_{n
\in \mathbb{N}}$ thus (\ref{Main result-ccl}) is proved. \endproof
\\

\textit{Proof of Corollary~\ref{Cor-sigma non nul} : } Let $C^{*},
\sigma^{*}>0$ be as in Proposition~\ref{Vap-analytic}. There exists
$N=N(\Gamma,s,\epsilon,\gamma) \in \NN^{*}$ large enough so that
\begin{equation} \label{Def-N}
\frac{\Gamma^{2}}{\left( \lambda_{N+1}-\frac{C^{*}\sigma^{*2}}{N+1} \right)^{s}}
\leqslant
\frac{\epsilon \gamma^{2}}{1-\epsilon}.
\end{equation}
Let $\sigma^{**}=\sigma^{**}(N)$ be as in Theorem~\ref{Thm-Main
result}. (notice that $\sigma^{**} \leqslant \sigma^{*}$) and
$\sigma \in (-\sigma^{**},\sigma^{**})-\{0\}$. Let $\Psi_{0} \in
H^{s}_{(\sigma)}(I,\CC) \cap \SSS$ verifying (\ref{Hyp-Cor-Main
result}). In order to get the conclusion of Corollary~\ref{Cor-sigma
non nul}, we prove that (\ref{Hyp-Main result}) holds, and we apply
the Theorem~\ref{Thm-Main result}. Using (\ref{vap-borne-0}), we get
$$\begin{array}{ll}
\sum\limits_{k=N+1}^{\infty} |\langle \Psi_{0} , \phi_{k,\sigma} \rangle |^{2}
& \leqslant
\frac{1}{\lambda_{N+1,\sigma}^{s}} \sum_{k=N+1}^{\infty}
\lambda_{k,\sigma}^{s} |\langle \Psi_{0} , \phi_{k,\sigma} \rangle |^{2}
\\ & \leqslant
\frac{1}{\lambda_{N+1,\sigma}^{s}} \sum_{k=1}^{\infty}
\lambda_{k,\sigma}^{s} |\langle \Psi_{0} , \phi_{k,\sigma} \rangle
|^{2}
\\ & \leqslant
\frac{\Gamma^{2}}{\left( \lambda_{N+1} - \frac{C^{*}\sigma^{2}}{N+1}
\right)^{s}}.
\end{array}$$
Thus (\ref{Def-N}) implies (\ref{Hyp-Main result}). \endproof

\section{Stabilization of $\CCC_{1}$}\label{sec:main}
\label{sigma=0} In all this section, the constants $C^{*},
\sigma^{*}$ are as in Proposition~\ref{Vap-analytic}.

\subsection{Main result}
\label{Main result-sigma=0} First, let us state the existence of an
implicit function $\sigma(\Psi)$ that will be used in the feedback
law. When $X$ is a normed space, $a \in X$ and $r>0$, we use the
notation $B_{X}(a,r):=\{ y \in X ; \|y-a\|_{X} <r \}$.
\begin{proposition} \label{Existence-sigma}
Let $N \in \NN^{*}$, $\epsilon>0$, and
$\theta \in C^{\infty}(\RR_{+},[0,\sigma^{*}])$ be such that
\begin{equation} \label{hyp-theta}
\begin{array}{lll}
\theta (0)=0 , & \theta(s)>0 ~~~\forall s>0 ,& \| \theta'
\|_{L^{\infty}} \leqslant \frac{1}{36NC^{*}}.
\end{array}
\end{equation}
There exists a unique
$\sigma \in C^{\infty}( B_{L^{2}}(0,2) , [0,\|\theta\|_{L^{\infty}}])$
such that
$$\sigma(\psi)=\theta(
\mathcal{V}_{\sigma(\psi),N,\epsilon}(\psi)),\qquad \forall \psi \in
B_{L^{2}}(0,2),$$ where $\mathcal{V}_{\sigma,N,\epsilon}$ is defined
by (\ref{Def-V-Lyap}).
\end{proposition}

The proof of this proposition is done in \cite{beauchard-et-al04}.
For sake of completeness, we repeat it in the Appendix. The main
result of this section is the following.

\begin{theorem} \label{Thm-Main result-0}
Let $N \in \NN^{*}$, $\gamma \in (0,1)$, $\epsilon>0$, $\theta \in
C^{\infty}(\RR_{+},[0,\sigma^{*}])$ verifying (\ref{hyp-theta}),
\begin{equation} \label{hyp-theta-BIS}
\|\theta\|_{L^{\infty}} \leqslant \min \left\{ \frac{1}{C^{*}}
\left( \frac{\epsilon \gamma^{2} N}{32(1-\epsilon/2)} \right)^{1/2},
\frac{\gamma}{2C^{*}}, \sigma^{\sharp}(N), \frac{1}{C^{*}} (
\sqrt{1-\epsilon/2} - \sqrt{1-\epsilon} ) \right\}
\end{equation}
and
\begin{equation} \label{Hyp-theta-Lem}
\| \theta' \|_{L^{\infty}} < \frac{1}{3(1+NC^{*})}.
\end{equation}
Let $\sigma \in C^{\infty}( B_{L^{2}}(0,2) ,
[0,\|\theta\|_{L^{\infty}}])$ be as in Proposition
\ref{Existence-sigma}. For every $\Psi_{0} \in \SSS \cap (H^{2}\cap
H^{1}_{0})(I,\CC)$ with
\begin{equation} \label{Hyp-psi0-0}
\begin{array}{lll}
\sum_{k=N+1}^{\infty} |\langle \Psi_{0} , \phi_{k} \rangle |^{2}
< \frac{\epsilon \gamma^{2}}{32(1-\epsilon/2)}
& \text{ and } &
|\langle \Psi_{0} , \phi_{1} \rangle | \geqslant \gamma,
\end{array}
\end{equation}
the Cauchy problem
(\ref{eq:main})-(\ref{eq:initial})-(\ref{eq:bords}). with
$u(t)=\sigma(\Psi(t))+v_{\sigma(\psi(t)),N,\epsilon}(\Psi(t))$ has a
unique strong solution $\psi$. Moreover this solution satisfies
\begin{equation} \label{CCL-lim}
\liminf\limits_{t \rightarrow + \infty}
| \langle \Psi(t) , \phi_{1} \rangle |^{2} \geqslant 1-\epsilon.
\end{equation}
\end{theorem}

The following Corollary states that the quantity $N$ appearing in
the feedback law may be uniform in a fixed bounded subset of $H^{s}$
for $s>0$.
\begin{corollary} \label{Cor-Main result-sigma=0}
Let $s>0$, $\epsilon>0$, $\Gamma>0$ and $\gamma \in (0,1)$. There
exists $N=N(\Gamma,s,\epsilon,\gamma) \in \NN^{*}$ such that, for
every $\Psi_{0} \in \SSS \cap (H^{2} \cap H^{1}_{0})(I,\CC)$ with $
\Psi_{0} \in H^{s}_{(0)}(I,\CC)$,
\begin{equation} \label{Hyp-Psi0-sigma=0}
\| \Psi_{0} \|_{H^{s}_{(0)}} \leqslant \Gamma
\text{  and  }
|\langle \Psi_{0} , \phi_{1} \rangle | \geqslant \gamma,
\end{equation}
the Cauchy problem
(\ref{eq:main})-(\ref{eq:initial})-(\ref{eq:bords}), with
$u(t)=\sigma(\Psi(t))+v_{\sigma(\psi(t)),N,\epsilon}(\Psi(t))$ has a
unique strong solution $\Psi$. Moreover this solution satisfies
(\ref{CCL-lim}).
\end{corollary}

\begin{remark}
The Theorem~\ref{Stab-Hs-sigma} with $\sigma=0$ is a direct
consequence of the Corollary~\ref{Cor-Main result-sigma=0}. The
feedback law, evoked in Theorem~\ref{Stab-Hs-sigma}, is implicitly
given by the Corollary~\ref{Cor-Main result-sigma=0}.
\end{remark}

\subsection{Preliminaries}
\label{Preliminaries-bis} \hspace{2 cm}\vspace{.5cm}

\begin{lemma} \label{Lem-diff-sigma}
Let $N \in \NN^{*}$, $\epsilon>0$ and $\theta$ satisfying
(\ref{hyp-theta}). There exist $C(N)>0$ and $\tilde{C}(N)>0$ such
that, for all $\xi_{1}, \xi_{2} \in B_{L^{2}}(0,1)$,
\begin{alignat}{4}
|\sigma(\xi_{1}) - \sigma(\xi_{2})| & \leqslant
3N\|\theta'\|_{L^{\infty}} \|\xi_{1}-\xi_{2}\|_{L^{2}},
\label{diff-sigma}\\
|\sigma(\xi_{1}) - \sigma(\xi_{2})| & \leqslant C(N)
\|\theta'\|_{L^{\infty}}
\|\xi_{1}-\xi_{2}\|_{H^{-1}},\label{diff-sigma-H^-1}\\
|v_{\sigma(\xi_{1}),N,\epsilon}(\xi_{1}) -
v_{\sigma(\xi_{2}),N,\epsilon}(\xi_{2}) | & \leqslant N(1+3NC^{*}
\|\theta'\|_{L^{\infty}}) \|\xi_{1}-\xi_{2}\|_{L^{2}},
\label{diff-v}\\
|v_{\sigma(\xi_{1}),N,\epsilon}(\xi_{1}) -
v_{\sigma(\xi_{2}),N,\epsilon}(\xi_{2}) | & \leqslant \tilde{C}(N)
\|\xi_{1}-\xi_{2}\|_{H^{-1}}. \label{diff-v-H1}
\end{alignat}
\end{lemma}

\textit{Proof of Lemma~\ref{Lem-diff-sigma} : } Since $N$ and
$\epsilon$ are fixed, in order to simplify the notations, we remove
them from the subscripts of this proof. We have
\begin{equation} \label{IAF}
|\sigma(\xi_{1}) - \sigma(\xi_{2})| \leqslant
\|\theta'\|_{L^{\infty}} |
\mathcal{V}_{\sigma(\xi_{1})}(\xi_{1}) - \mathcal{V}_{\sigma(\xi_{2})}(\xi_{2}) |.
\end{equation}
Using
\begin{equation} \label{decomposition}
\begin{array}{ll}
|\langle \xi_{1} ,\phi_{k,\sigma_{1}} \rangle |^{2} - |\langle
\xi_{2} ,\phi_{k,\sigma_{2}} \rangle |^{2} =& \langle   \xi_{1} -
\xi_{2}   ,   \phi_{k,\sigma_{1}}   \rangle \overline{ \langle
\xi_{1}   ,   \phi_{k,\sigma_{1}}   \rangle } \\
&
+\langle   \xi_{2} , \phi_{k,\sigma_{1}}   \rangle \overline{
\langle   \xi_{1} - \xi_{2} ,   \phi_{k,\sigma_{1}} \rangle }
\\ &
+\langle   \xi_{2}   ,   \phi_{k,\sigma_{1}} - \phi_{k,\sigma_{2}}   \rangle
\overline{ \langle   \xi_{2}   ,   \phi_{k,\sigma_{1}}  \rangle }
\\ &
+\langle   \xi_{2}   ,   \phi_{k,\sigma_{2}}   \rangle
\overline{ \langle \xi_{2} ,  \phi_{k,\sigma_{1}} - \phi_{k,\sigma_{2}}   \rangle }
\end{array}
\end{equation}
and (\ref{diff-vep-k}), we get
$$| \mathcal{V}_{\sigma(\xi_{1})}(\xi_{1}) - \mathcal{V}_{\sigma(\xi_{2})}(\xi_{2}) |
\leqslant
2N\|\xi_{1}-\xi_{2}\|_{L^{2}} +
2NC^{*}|\sigma(\xi_{1})-\sigma(\xi_{2})|,$$
$$| \mathcal{V}_{\sigma(\xi_{1})}(\xi_{1}) - \mathcal{V}_{\sigma(\xi_{2})}(\xi_{2}) |
\leqslant 2N C_{1}(N) \|\xi_{1}-\xi_{2}\|_{H^{-1}} +
2NC^{*}|\sigma(\xi_{1})-\sigma(\xi_{2})|.$$ where $C_{1}(N):= \max
\{ \| \varphi_{k,\sigma} \|_{H^{1}_{0}}  ; k \in \{1,\cdots,N\},
\sigma \in [0,\sigma^{*}] \}.$ Using the previous inequalities and
(\ref{hyp-theta}), we get
$$\frac{17}{18} | \sigma(\xi_{1}) - \sigma(\xi_{2}) | \leqslant
2N \|\theta'\|_{\infty} \|\xi_{1}-\xi_{2}\|_{L^{2}},$$
$$\frac{17}{18} | \sigma_{1} - \sigma_{2} | \leqslant
2N C_{1}(N) \mathcal{C} \|\theta'\|_{\infty}
\|\xi_{1}-\xi_{2}\|_{H^{-1}},
$$
which implies (\ref{diff-sigma}) and
(\ref{diff-sigma-H^-1}) with $C(N)=3N\mathcal{C}C_{1}(N)$.

Let us write $v_{j}$ instead of $v_{\sigma(\xi_{j})}(\xi_{j})$.
Using, for the term
$$\langle x \xi_{1} , \phi_{j,\sigma(\xi_{1})} \rangle
\overline{\langle \xi_{1} , \phi_{j,\sigma(\xi_{1})} \rangle} -
\langle x \xi_{2} , \phi_{j,\sigma(\xi_{2})} \rangle
\overline{\langle \xi_{2} , \phi_{j,\sigma(\xi_{2})} \rangle}$$ the
same kind of decomposition as in (\ref{decomposition}), together
with (\ref{diff-vep-k}), we get
$$|v_{1}-v_{2}| \leqslant
N \|\xi_{1}-\xi_{2}\|_{L^{2}} + NC^{*} |\sigma(\xi_{1})-\sigma(\xi_{2})|,$$
$$|v_{1}-v_{2}| \leqslant
2N\mathcal{C}C_{1}(N)\|\xi_{1}-\xi_{2}\|_{H^{-1}} + 2NC^{*}
|\sigma(\xi_{1})-\sigma(\xi_{2})|,$$ where $\mathcal{C}$ is defined
by (\ref{Def-cst-x}). Thus, using (\ref{diff-sigma}) and
(\ref{diff-sigma-H^-1}), we get (\ref{diff-v}) and (\ref{diff-v-H1})
with $\tilde{C}(N):= 2N[ \mathcal{C}C_{1}(N) + C^{*} C(N)
\|\theta'\|_{\infty}]$. \endproof

\begin{proposition} \label{Existence-0}
Let $N \in \NN^{*}$, $\epsilon>0$, $\theta$ verifying
(\ref{hyp-theta}) and~\eqref{Hyp-theta-Lem}. For every $\Psi_{0} \in
\SSS$ the Cauchy problem
(\ref{eq:main})-(\ref{eq:initial})-(\ref{eq:bords}) with
$u(t)=\sigma(\Psi(t))+v_{\sigma(\psi(t)),N,\epsilon}(\Psi(t))$ has a
unique weak solution i.e. $\Psi \in C^{0}(\mathbb{R},\SSS) \cap
C^{1}((0,+\infty),$ $H^{-2}_{(0)})$. If, moreover $\Psi \in (H^{2}
\cap H^{1}_{0})(I,\CC)$, then $\Psi$ is a strong solution i.e. $\Psi
\in C^{0}(\mathbb{R},H^{2} \cap H^{1}_{0}) \cap
C^{1}((0,+\infty),L^{2})$.
\end{proposition}

\textit{Proof of Proposition~\ref{Existence-0} : } The strategy is
the same as in the proof of Proposition~\ref{Existence}. Let $T>0$
be such that
$$NTe^{(N+\|\theta\|_{L^{\infty}})T} < \frac{1}{2}.$$
Let $\Psi_{0} \in \SSS$. In order to build solutions on $[0,T]$,
we apply the Banach fixed point theorem to the map
$$\begin{array}{cccc}
\Theta : & C^{0}([0,T],\SSS)    & \rightarrow & C^{0}([0,T],\SSS) \\
         & \xi                  & \mapsto     & \Psi
\end{array}$$
where $\Psi$ is the weak solution of
(\ref{eq:main})-(\ref{eq:initial})-(\ref{eq:bords}) with
$u(t)=\sigma(\xi(t)) + v_{\sigma(\xi(t)),N,\epsilon}(\xi(t))$.

The map $\Theta$ is well defined and maps $C^{0}([0,T],\SSS)$ into
itself, moreover, it takes values in $C^{0}([0,T],\SSS) \cap
C^{1}((0,T),H^{-2}_{(0)})$ (see
Proposition~\ref{Existence-classic}). Let us prove that $\Theta$ is
a contraction of $C^{0}([0,T],\SSS)$. Let $\xi_{j} \in
C^{0}([0,T],\SSS)$,
$v_{j}:=v_{\sigma(\xi_{j}),N,\epsilon}(\xi_{j})$,
$\Psi_{j}:=\Theta(\xi_{j})$, for $j=1,2$ and
$\Delta:=\Psi_{1}-\Psi_{2}$. We have
$$
\Delta(t)=i \int_{0}^{t}
e^{-iA(t-s)} [
(\sigma(\xi_{1})+v_{1}) x \Delta (s) +
(\sigma(\xi_{1}) - \sigma(\xi_{2}) + v_{1}-v_{2}) x \Psi_{2}(s) ]ds.
$$
Using (\ref{diff-sigma}) and (\ref{diff-v}), we get
$$\begin{array}{ll}
\|\Delta(t)\|_{L^{2}} \leqslant
&
\int_{0}^{t}
\Big( \|\theta'\|_{L^{\infty}} + N \Big) \|\Delta(s)\|_{L^{2}} ds
\\ &
+ \int_{0}^{t}
\Big( 3 N \|\theta'\|_{L^{\infty}} + N[ 1+3NC^{*} \|\theta'\|_{L^{\infty}} ] \Big)
\| \xi_{1}-\xi_{2} \|_{L^{2}} ds.
\end{array}$$
Thus, the Gronwall Lemma implies
$$
\|\Delta\|_{C^{0}([0,T],L^{2})} \leqslant
\|\xi_{1}-\xi_{2}\|_{C^{0}([0,T],L^{2})}
[ 1 + 3 (1+NC^{*}) \|\theta'\|_{L^{\infty}} ] NT e^{T[N+\|\theta\|_{L^{\infty}}]}.
$$
The choice of $T$ and (\ref{Hyp-theta-Lem}) ensure that $\Theta$ is
a contraction of $C^{0}([0,T],\SSS)$. Therefore, there exists a
fixed point $\Psi \in C^{0}([0,T],\SSS)$ such that
$\Theta(\Psi)=\Psi$. Since $\Theta$ takes values in
$C^{0}([0,T],\SSS) \cap C^{1}([0,T],H^{-2}_{(0)})$, necessarily
$\Psi$ belongs to this space, thus, it is a weak solution of
(\ref{eq:main})-(\ref{eq:initial})-(\ref{eq:bords}) on $[0,T]$.

If, moreover, $\Psi_{0} \in (H^{2} \cap H^{1}_{0})(I,\CCC)$, then
the map $\Theta$ takes values in $C^{0}([0,T],H^{2} \cap H^{1}_{0})
\cap C^{1}([0,T],L^{2})$ thus $\Psi$ belongs to this space and it is
a strong solution.

Since the time $T$ does not depend on $\Psi_{0}$, the solution can
be continued globally in time. We, therefore, have the existence of
global solutions to the closed-loop system. \endproof

\begin{proposition} \label{Convergence-H^{-1}-0}
Let $\sigma>0$, $N \in \NN$, $\epsilon >0$, $\theta$ as
in~\eqref{hyp-theta}, $(\Psi_{0}^{n})_{n \in \NN}$ be a sequence of
$\SSS$ and $\Psi_{0}^{\infty} \in L^{2}$ with
$\|\Psi^{\infty}_0\|_{L^{2}} \leqslant 1$ be such that
$$\lim\limits_{n \rightarrow + \infty} \Psi_{0}^{n} = \Psi_{0}^{\infty}
\text{ strongly in } H^{-1}(I,\CC).$$
Let $\Psi^{n}$ (resp. $\Psi^{\infty}$) be the weak solution of
(\ref{eq:main})-(\ref{eq:initial})-(\ref{eq:bords})
with $u(t)=\sigma(\Psi^{n}(t))+v_{\sigma(\Psi^{n}(t)),N,\epsilon}(\Psi^{n}(t))$
(resp. with  $u(t)=\sigma(\Psi^{\infty})+
v_{\sigma(\Psi^{\infty}),N,\epsilon}(\Psi^{\infty}(t))$).
Then, for every $\tau>0$,
$$\lim\limits_{n \rightarrow + \infty} \Psi^{n}(\tau) = \Psi^{\infty}(\tau)
\text{ strongly in } H^{-1}(I,\CC).$$
\end{proposition}

\textit{Proof of Proposition~\ref{Convergence-H^{-1}-0} :} The proof
exactly follows that of the Proposition~\ref{Convergence-H^{-1}}. In
order to simplify the notations, we write $v(\Psi)$ instead of
$v_{\sigma(\Psi),N,\epsilon}(\Psi)$. We have
$$\begin{array}{ll}
(\Psi^{n}-\Psi^{\infty})(t)=
&
e^{-iAt} ( \Psi^{n}_{0} - \Psi^{\infty}_{0} )
+i \int_{0}^{t} e^{-iA(t-s)}
[ \sigma(\Psi^{n})-\sigma(\Psi^{\infty}) ]  x \Psi^{n} ds
\\ &
+i \int_{0}^{t} e^{-iA(t-s)}
[ v(\Psi^{n})-v(\Psi^{\infty}) ] x \Psi^{n} ds
\\ &
+i \int_{0}^{t} e^{-iA(t-s)}
[ \sigma(\Psi^{\infty}) + v(\Psi^{\infty}) ] x (\Psi^{n}-\Psi^{\infty}) ds.
\end{array}$$
Using (\ref{diff-sigma-H^-1}), (\ref{diff-v-H1}) and
$\|x\Psi\|_{H^{-1}} \leqslant \|x\Psi\|_{L^{2}} \leqslant 1$, we get
$$\begin{array}{ll}
&\| (\Psi^{n}-\Psi^{\infty})(t) \|_{H^{-1}}
 \leqslant
\| \Psi^{n}_{0} - \Psi^{\infty}_{0} \|_{H^{-1}}
\\&
+ \int_{0}^{t} \Big(
C(N)\|\theta'\|_{L^{\infty}} + \tilde{C}(N)
+ \mathcal{C}(\|\theta\|_{L^{\infty}}+N)
\Big)
\| \Psi^{n}-\Psi^{\infty} \|_{H^{-1}} ds,
\end{array}$$
where $\mathcal{C}$ is given by (\ref{Def-cst-x}). The Gronwall
lemma concludes the proof. \endproof

\subsection{Proof of Theorem~\ref{Thm-Main result-0} and Corollary
\ref{Cor-Main result-sigma=0}} \label{subsection:Cor-Main
result-sigma=0} \hspace{2 cm}\vspace{.5 cm}

\textit{Proof of Theorem~\ref{Thm-Main result-0} : } For $\varphi
\in B_{L^{2}}(0,2)$, we define
$$\mathcal{V}_{N,\epsilon}(\varphi):=
\mathcal{V}_{\sigma(\varphi),N,\epsilon}(\varphi),$$
where $\mathcal{V}_{\sigma,N,\epsilon}$ is defined by (\ref{Def-V-Lyap}).
Since $N$ and $\epsilon$ are fixed, in order to simplify the notations,
we omit them in the subscripts of this proof, and we write
$v(\Psi)$ instead of $v_{\sigma(\Psi),N,\epsilon}(\Psi)$.

Let $\Psi_{0} \in \SSS \cap (H^{2} \cap H^{1}_{0})(I,\CC)$ and
$\Psi$ be the strong solution of
(\ref{eq:main})-(\ref{eq:initial})-(\ref{eq:bords}) with
$u(t)=\sigma(\Psi(t))+v_{\sigma(\psi(t)),N,\epsilon}(\Psi(t))$ given
by Proposition~\ref{Existence-0}. Since $\Psi \in
C^{1}(\mathbb{R},L^{2})$ and $\sigma \in
C^{\infty}(B_{L^{2}}(0,2))$, the map $t \mapsto
\mathcal{V}(\psi(t))$ is $C^{1}$. We have
$$\frac{d}{dt} \mathcal{V}(\Psi) =
- 2 v(\Psi)^{2} - \frac{d}{dt} \Big[ \sigma(\Psi) \Big] \Re \left(
\sum_{k=1}^{N} a_{k} \langle \Psi , \frac{d\phi_{k,\sigma}}{d\sigma}
\Big|_{\sigma(\Psi)} \rangle \overline{\langle \Psi ,
\phi_{k,\sigma(\Psi)} \rangle} \right),$$ where $a_{1}:=1$ and
$a_{k}:=1-\epsilon$ for $k=2,\cdots,N$. Moreover,
$$\frac{d}{dt} \Big[ \sigma(\Psi) \Big]
=
\theta'(\mathcal{V}(\psi)) \frac{d}{dt} \mathcal{V}(\Psi)$$
thus
\begin{equation} \label{dV/dt-sigma=0}
\left[
1+2 \theta'(\mathcal{V}(\psi))
\Re \left( \sum_{k=1}^{N} a_{k}
\langle \Psi , \frac{d\phi_{k,\sigma}}{d\sigma} \Big|_{\sigma(\Psi)} \rangle
\overline{\langle \Psi , \phi_{k,\sigma(\Psi)} \rangle} \right)
\right]
\frac{d}{dt} \mathcal{V}(\Psi)
= - 2 v(\Psi)^{2}.
\end{equation}
Using (\ref{vep-k-1}) and (\ref{hyp-theta}), we get
$$1+2 \theta'(\mathcal{V}(\psi))
\Re \left( \sum_{k=1}^{N} a_{k}
\langle \Psi , \frac{d\phi_{k,\sigma}}{d\sigma} \Big|_{\sigma(\Psi)} \rangle
\overline{\langle \Psi , \phi_{k,\sigma(\Psi)} \rangle} \right)
\geqslant
1-2\|\theta'\|_{L^{\infty}}NC^{*}
>0$$
thus, $t \mapsto \mathcal{V}(\Psi(t))$ is a non increasing function.
There exists $\alpha \in [0,\mathcal{V}(\Psi_{0})]$ such that
$$\lim\limits_{t \rightarrow + \infty}
\mathcal{V}(\Psi(t)) = \alpha.$$ Using (\ref{diff-vep-k}),
(\ref{hyp-theta-BIS}) and (\ref{Hyp-psi0-0}), we get
$$\begin{array}{ll}
|\langle \Psi_{0} , \phi_{1,\sigma(\Psi_{0})} \rangle|
&
\geqslant
|\langle \Psi_{0} , \phi_{1} \rangle|
- |\langle \Psi_{0} , \phi_{1} - \phi_{1,\sigma(\Psi_{0})} \rangle|
\\  & \geqslant
\gamma - C^{*} \|\theta\|_{\infty}
\\ & \geqslant \tilde{\gamma} := \frac{\gamma}{2},
\end{array}$$
$$\begin{array}{ll}
\sum_{k=N+1}^{\infty}  |\langle \Psi_{0} , \phi_{k,\sigma(\Psi_{0})}
\rangle |^{2} & \leqslant 2 \sum_{k=N+1}^{\infty} \left(|\langle
\Psi_{0} , \phi_{k} \rangle |^{2} + |\langle \Psi_{0} ,
\phi_{k,\sigma(\Psi_{0})} - \phi_{k} \rangle |^{2}\right)
\\ & \leqslant
\frac{\epsilon \gamma^{2}}{16(1-\epsilon/2)} +
2(C^{*} \|\theta\|_{L^{\infty}})^{2} \sum_{k=N+1}^{\infty} \frac{1}{k^{2}}
\\ & \leqslant
\frac{\epsilon \gamma^{2}}{16(1-\epsilon/2)} +
\frac{2(C^{*}\|\theta\|_{L^{\infty}})^{2}}{N}
\\ & \leqslant
\frac{\tilde{\epsilon} \tilde{\gamma}^{2}}{(1-\tilde{\epsilon})}
\end{array}$$
where $\tilde{\epsilon}:=\epsilon/2$. Thus, as in the proof of
Theorem ~\ref{Thm-Main result},
$\mathcal{V}(\Psi_{0})<\tilde{\epsilon}$, so $\alpha \in
(0,\tilde{\epsilon})$.

Let $(t_{n})_{n \in \mathbb{N}}$ be an increasing sequence
of positive real numbers such that $t_{n} \rightarrow + \infty$
when $n \rightarrow + \infty$. Since $\|\Psi(t_{n})\|_{L^{2}}=1$
for every $n \in \mathbb{N}$, there exists
$\Psi_{\infty} \in L^{2}(I,\CC)$ such that, up to an extraction
$$\Psi(t_{n}) \rightarrow \Psi_{\infty}
\text{ weakly in } L^{2}(I,\CC)
\text{ and strongly in } H^{-1}(I,\CC).$$
Let $\xi$ be the weak solution of
$$\left\lbrace \begin{array}{l}
i \frac{\partial \xi}{\partial t} = A_{\sigma} \xi
- v_{\sigma(\xi),N,\epsilon}(\xi(t)) x \xi, \\
\xi(t,\pm 1/2)=0,\\
\xi(0)=\Psi_{\infty}.
\end{array} \right.$$
Thanks to the Proposition~\ref{Convergence-H^{-1}-0}, for every
$\tau>0$, $\Psi(t_{n}+\tau) \rightarrow \xi(\tau)$ strongly in
$H^{-1}(I,\mathbb{C})$ when $n \rightarrow + \infty$, thus
$\sigma(\Psi(t_{n}+\tau)) \rightarrow \sigma(\xi(\tau))$ when $n
\rightarrow + \infty$ (see Lemma~\ref{Lem-diff-sigma}). Therefore,
$\mathcal{V}(\Psi(t_{n}+\tau)) \rightarrow \mathcal{V}(\xi(\tau))$
when $n \rightarrow + \infty$, so $\mathcal{V}(\xi) \equiv \alpha$.
Thus, $\sigma(\xi) \equiv \overline{\sigma}:=\theta(\alpha)$ and we
have, for every $t \in \mathbb{R}_{+}$,
$$\mathcal{V}(\xi(t))=1-|\langle \xi(t) , \phi_{1,\overline{\sigma}} |^{2}
- (1-\epsilon) \sum_{k=2}^{N} |\langle \xi(t) , \phi_{k,\overline{\sigma}} \rangle |^{2}.$$
Since $\xi \in C^{1}(\mathbb{R}_{+},H^{-2}_{(0)})$, the previous equality implies
$$\frac{d\mathcal{V}(\xi)}{dt} = -2v(\xi)^{2}.$$
Since $\mathcal{V}(\xi) \equiv \alpha$, then $v(\xi) \equiv 0$.

\emph{First case : $\alpha=0$.} Then $\mathcal{V}(\Psi(t))
\rightarrow 0$ when $t \rightarrow + \infty$ and
$\overline{\sigma}=0$. Moreover, for every $t \in (0,\infty)$,
$$\begin{array}{ll}
\VVV(\Psi(t)) & \geqslant 1 - |\langle \Psi , \phi_{1,\sigma(\Psi)}
\rangle |^{2} -(1-\epsilon) \sum_{k=2}^{\infty} |\langle \Psi ,
\phi_{k,\sigma(\Psi)} \rangle |^{2}
\\ &
\geqslant
\epsilon ( 1 - |\langle \Psi , \phi_{1,\sigma(\Psi)} \rangle |^{2} ),
\end{array}$$
Thus,
$$|\langle \Psi(t_{n}) , \phi_{1,\sigma(\Psi(t_{n}))} \rangle |
\rightarrow 1,$$
which leads to
$$|\langle \Psi(t_{n}) , \phi_{1} \rangle |
\rightarrow 1.$$ because $\sigma(\Psi(t_{n})) \rightarrow 0$.

\emph{Second case : $\alpha \neq 0$.} Then
$\overline{\sigma}=\theta(\alpha) >0$. Exactly as in the first
analysis, done in the proof of Theorem~\ref{Thm-Main result}, we get
$$\Psi_{\infty} = \beta \phi_{1,\overline{\sigma}}$$
where $\beta \in \mathbb{C}$ and $|\beta|^{2}>1-\tilde\epsilon$.
Thus
$$
\lim\limits_{n \rightarrow + \infty} | \langle \Psi(t_{n}) ,
\phi_{1} \rangle |  = |\langle \Psi_{\infty} , \phi_{1} \rangle |
\geqslant |\beta| - |\langle \Psi_{\infty} ,
\phi_{1,\overline{\sigma}} - \phi_{1} \rangle | \geqslant
\sqrt{1-\epsilon/2} - C^{*} \overline{\sigma}
$$
where we used (\ref{vap-borne-0}) in the last inequality. Finally,
thanks to $0<\overline{\sigma} \leqslant \|\theta\|_{\infty}$ and
(\ref{hyp-theta-BIS}), we get (\ref{CCL-lim}). \endproof
\\

\textit{Proof of Corollary~\ref{Cor-Main result-sigma=0} :} It can
be done in a very similar way to the proof of the Corollary
\ref{Cor-sigma non nul}. \endproof

\section{Numerical simulations}\label{sec:numeric}
\begin{figure}[t]
\psfrag{L}{\tiny{Lyapunov
function}}\psfrag{V}{$\VVV_\epsilon(\Psi)$}
\psfrag{D}{$1-|\her{\Psi}{\phi_{1,\sigma}}|^2$}\psfrag{C}{\tiny{Feedback
law}}\psfrag{U}{$v_\epsilon(\Psi)$}\psfrag{T}{\tiny{Time}}
\begin{center}
\includegraphics[height=6.5 cm, width=9 cm]{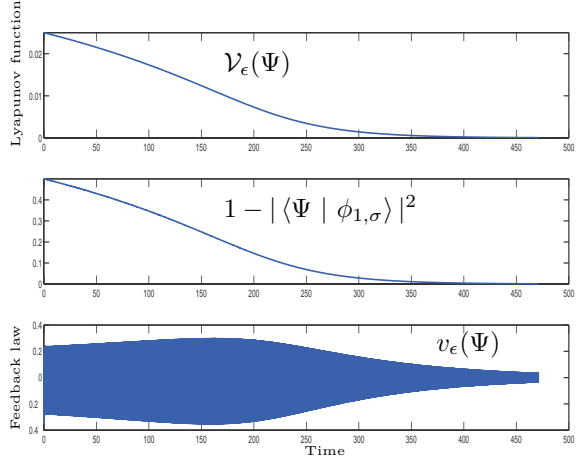}
\caption{The approximate stabilization of $\CCC_{1,\sigma}$, where
$\Psi_0=\frac{1}{\sqrt{2}}(\phi_{1,\sigma}+\phi_{3,\sigma})$ and
therefore the cut-off dimension is 3; as it can be seen the
closed-loop system reaches the .05-neighborhood of $\phi_{1,\sigma}$
in a time $T=150\pi$ corresponding to about 200 periods of the
longest natural period corresponding to the ground to the first
excited state.}
\label{fig1}                                 
\end{center}                                 
\end{figure}
In this section, we check out the performance of the techniques on
some numerical simulations. We consider, as a test case, the
stabilization of the initial state
$\Psi_0=\frac{1}{\sqrt{2}}(\phi_{1,\sigma}+\phi_{3,\sigma})$ around
the ground state $\phi_{1,\sigma}$. Therefore, the cut-off dimension
$N$ is $3$. Note that, such a test case is particularly a hard one
in a near-degenerate situation. Indeed, considering the feedback
law~\eqref{eq:law1} for $\sigma=0$, one can easily see that for
parity reasons $v(\Psi(t))\equiv 0$.

In a first simulation, we consider the non-degenerate case of
$\sigma\neq 0$. As mentioned above the constant $\sigma$ needs to be
small. In fact, one should choose $\sigma$, such that the
perturbation $\sigma x$ is small compared to the the operator
$-\frac{1}{2}\frac{\partial^2}{\partial x^2}$. We choose it here to
be $\sigma=2e+01$. The Figure~\ref{fig1} illustrates the simulation
of the closed-loop system when $u=\sigma+v_\epsilon$ with
$\varsigma=1e+03$ and $\epsilon=5e-02$. The simulations have been
done applying a third order split-operator method (see
e.g.~\cite{Feit-et-al-82}), where instead of computing
$\exp(-i~dt~(A_\sigma-v_\epsilon x))$ at each time step, we compute
$$
\exp(-i~dt~A_\sigma/2)\exp(i~dt~v_\epsilon x)\exp(-i~dt~A_\sigma/2).
$$
Moreover, we consider a Galerkin discretization over the first 20
modes of the system (it turns out, by considering higher modal
approximations, that 20 modes are completely sufficient to get a
trustable result).

\begin{figure}[t]
\psfrag{L}{\tiny{Lyapunov
function}}\psfrag{V}{$\VVV_\epsilon(\Psi)$}
\psfrag{D}{$1-|\her{\Psi}{\phi_{1}}|^2$}\psfrag{C}{\tiny{Feedback
law}}\psfrag{U}{$\sigma(\Psi)+v_\epsilon(\Psi)$}\psfrag{T}{\tiny{Time}}
\begin{center}
\includegraphics[height=6.5 cm, width=9 cm]{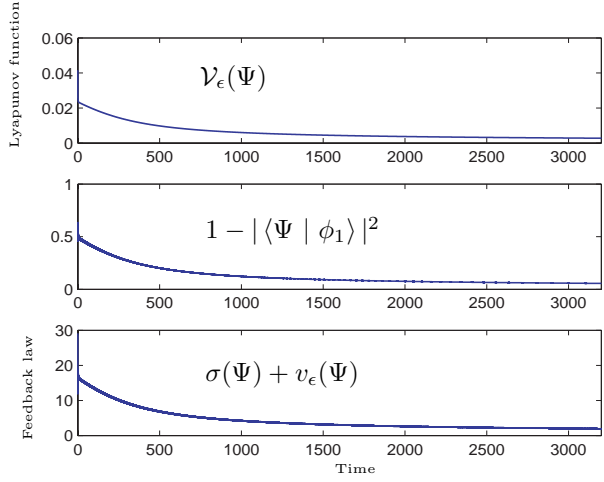}
\caption{The approximate stabilization of $\CCC_{1}$, where
$\Psi_0=\frac{1}{\sqrt{2}}(\phi_{1}+\phi_{3})$ and therefore the
cut-off dimension is 3; as it can be seen, the closed-loop system
reaches the .05-neighborhood of $\phi_{1}$ in a time $T=1000\pi$
corresponding to about 1300 periods of the longest natural period
corresponding to the ground to the first excited state.}
\label{fig2}                                 
\end{center}                                 
\end{figure}

Now, let us consider the degenerate case of $\sigma=0$. As mentioned
above, such a case is not treatable with the explicit feedback
design of~\eqref{eq:law1}. However, the simulations of
Figure~\ref{fig2}, show that the implicit Lyapunov design provided
in Subsection~\ref{ssec:heur} removes the degeneracy problem and
ensures the approximate stabilization of the initial state
$\frac{1}{\sqrt 2}(\phi_1+\phi_3)$ around the ground state $\phi_1$.

We consider the function $\theta(r)=\eta r$ with $\eta=7e+02$.
Furthermore, in the feedback design $v_\epsilon$, we consider
$\varsigma=1e+03$ and $\epsilon=5e-02$. The numerical scheme is
similar to the simulations of Figure~\ref{fig1}. In order to
calculate the implicit part of the feedback design $\sigma(\Psi)$,
we apply a fixed point algorithm.

\section{Appendix}

This appendix is devoted to the proofs of the
Proposition~\ref{Existence-classic} and the Proposition
\ref{Existence-sigma}.
\subsection{Proof of Proposition~\ref{Existence-classic} }
Let $\Psi_{0} \in \SSS$, $T_{1}>0$ and $u \in C^{0}([0,T_{1}],\RR)$.
Let $T \in (0,T_{1})$ be such that
\begin{equation} \label{u-contract}
\|u\|_{L^{1}(0,T)}<1.
\end{equation}
We prove the existence of
$\Psi \in C^{0}([0,T],L^{2}(I,\mathbb{C}))$ such that
(\ref{weak}) holds by applying the Banach fixed point theorem
to the map
$$\begin{array}{cccc}
\Theta : & C^{0}([0,T],L^{2}) & \rightarrow & C^{0}([0,T],L^{2})
\\
            &   \xi              & \mapsto     & \Psi
\end{array}$$
where $\Psi$ is the weak solution of
$$\left\lbrace\begin{array}{l}
i \frac{\partial \Psi}{\partial t} = A \Psi - u(t)x \xi , \\
\Psi(0,x)=\Psi_{0}(x),\\
\Psi(t,\pm 1/2)=0.
\end{array}\right.$$
i.e. $\Psi \in C^{0}([0,T],L^{2})$ and satisfies, for every $t \in [0,T]$,
$$\Psi(t)=e^{-iAt} \Psi_{0} +i\int_{0}^{t} e^{-iA(t-s)} u(s) x \xi(s) ds
\text{  in  } L^{2}(I,\CC).$$
Notice that $\Theta$ takes values in $C^{1}([0,T],H^{-2}_{(0)}(I,\CC))$.

For $\xi_{1}, \xi_{2} \in C^{0}([0,T],L^{2}(I,\CC))$,
$\Psi_{1}:=\Theta(\xi_{1})$,
$\Psi_{2}:=\Theta(\xi_{2})$
we have
$$(\Psi_{1}-\Psi_{2})(t)=i\int_{0}^{t} e^{-iA(t-s)} u(s) x(\xi_{1}-\xi_{2})(s) ds$$
thus
$$\|(\Psi_{1}-\Psi_{2})(t)\|_{L^{2}} \leqslant \int_{0}^{t} |u(s)| ds
\|\xi_{1}-\xi_{2}\|_{C^{0}([0,T],L^{2})}.$$
The assumption (\ref{u-contract}) guarantees that $\Theta$ is a contraction
of $C^{0}([0,T],L^{2})$, thus, $\Theta$ has a fixed point
$\Psi \in C^{0}([0,T],L^{2})$. Since $\Theta$ takes values in $C^{1}([0,T],H^{-2}_{(0)})$,
then $\Psi$ belongs to this space. Moreover, this function satisfies (\ref{weak}).

Finally, we have built weak solutions on $[0,T]$ for every $\Psi_{0}$,
and the time $T$ does not depend on $\Psi_{0}$, thus,
this gives solutions on $[0,T_{1}]$.

Let us prove that this solution is continuous with respect to the
the initial condition $\Psi_{0}$, for the $L^{2}(I,\CC)$-topology.
Let $\Psi_{0}, \Phi_{0} \in \SSS$ and $\Psi$, $\Phi$ the associated
weak solutions. We have
$$\|(\Psi-\Phi)(t)\|_{L^{2}} \leqslant
\| \Psi_{0} - \Phi_{0} \|_{L^{2}} +
\int_{0}^{t} |u(s)| \|(\Psi-\Phi)(s)\|_{L^{2}} ds,$$
thus Gronwall Lemma gives
$$\|(\Psi-\Phi)(t)\|_{L^{2}} \leqslant
\| \Psi_{0} - \Phi_{0} \|_{L^{2}} e^{\|u\|_{L^{1}(0,T_{1})}}.$$ This
gives the continuity of the weak solutions with respect to the
initial conditions.
\\

Now, let us assume that $\Psi_{0} \in H^{2} \cap H^{1}_{0}(I,\CC)$.
Take $C$ to be a positive constant such that for every $\varphi \in
H^{2}\cap H^{1}_{0}(I,\CC)$, $\|x\varphi \|_{H^{2}\cap H^{1}_{0}}$
$\leqslant C \|\varphi\|_{H^{2}\cap H^{1}_{0}}$. We consider, then,
$T>0$ such that  $C \|u\|_{L^{1}(0,T)}$ $< 1$. By applying the fixed
point theorem on
$$\Theta_{2}:C^{0}([0,T],H^{2}\cap H^{1}_{0}) \rightarrow
C^{0}([0,T],H^{2}\cap H^{1}_{0})$$ defined by the same expression as
$\Theta$, and using the uniqueness of the fixed point of $\Theta$,
we get that the weak solution is a strong solution. The continuity
with respect to the initial condition of the strong solution can
also be proved applying the same arguments as in above.
\\

Finally, let us justify that the weak solutions take values in
$\SSS$. For $\Psi_{0} \in H^{2}\cap H^{1}_{0}$, the solution belongs
to $C^{1}([0,T],L^{2}) \cap C^{0}([0,T],H^{2}\cap H^{1}_{0})$ thus,
the following computations are justified
$$\frac{d}{dt} \|\Psi(t)\|_{L^{2}}^{2} = 2 \Re
\langle \frac{\partial \Psi}{\partial t} , \Psi \rangle
=0.$$
Thus $\Psi(t) \in \SSS$ for every $t\in [0,T]$.

For $\Psi_{0} \in \SSS$, we get the same conclusion thanks to a
density argument and the continuity for the
$C^{0}([0,T],L^{2})$-topology of the weak solutions with respect to
the initial condition.
\endproof

\subsection{Proof of Proposition~\ref{Existence-sigma} }
Let $\Psi \in B_{L^{2}}(0,2)$. We prove the existence of
$\sigma(\Psi)$ by applying the Banach fixed point Theorem to the map
$$\begin{array}{cccc}
\Pi : & [0,\|\theta\|_{L^{\infty}}] & \rightarrow & [0,\|\theta\|_{L^{\infty}}] \\
      &  \sigma        & \mapsto     & \theta( \mathcal{V}_{\sigma,N,\epsilon}(\Psi))
\end{array}$$
For $\sigma_{1}, \sigma_{2} \in [0,\|\theta\|_{L^{\infty}}]$, we have
$$|\Pi(\sigma_{1}) - \Pi(\sigma_{2}) | \leqslant
\|\theta'\|_{L^{\infty}} | \mathcal{V}_{\sigma_{1},N,\epsilon}(\Psi)
- \mathcal{V}_{\sigma_{2},N,\epsilon}(\Psi)|.$$ Using the following
inequality\small
\begin{align*}
\Big| | \langle \Psi , \phi_{j,\sigma_{1}} \rangle |^{2}
- | \langle \Psi , \phi_{j,\sigma_{2}} \rangle |^{2} \Big|
 \leqslant &
\Big| \langle \Psi , \phi_{j,\sigma_{1}} - \phi_{j,\sigma_{2}}
\rangle \overline{  \langle \Psi , \phi_{j,\sigma_{1}} \rangle  }
\Big| +\Big| \langle \Psi , \phi_{j,\sigma_{2}} \rangle \overline{
\langle \Psi , \phi_{j,\sigma_{1}} - \phi_{j,\sigma_{2}} \rangle }
\Big|
\\ \leqslant &
8 \|\phi_{j,\sigma_{1}} - \phi_{j,\sigma_{2}} \|_{L^{2}},
\end{align*}\normalsize
together with (\ref{diff-vep-k}), we get
$$|\Pi(\sigma_{1}) - \Pi(\sigma_{2}) | \leqslant
8NC^{*}\|\theta'\|_{L^{\infty}} |\sigma_{1} - \sigma_{2}|.$$ Thus,
the assumption (\ref{hyp-theta}) ensures that $\Pi$ is a contraction
of $[0,\|\theta\|_{L^{\infty}}]$. Therefore, $\Pi$ has a unique
fixed point $\sigma(\Psi)$.

Now, let us prove that $\sigma$ is $C^{\infty}$. The map
$$\begin{array}{cccc}
F: & [0,\|\theta\|_{L^{\infty}}] \times B_{L^{2}}(0,2) & \rightarrow & \RR \\
   &                       (\sigma , \Psi)             & \mapsto     &
\sigma - \theta( \mathcal{V}_{\sigma,N,\epsilon}(\Psi))
\end{array}$$
is regular with respect to $\sigma$ and $\Psi$,
$F(\sigma(\Psi),\Psi)=0$, for every $\Psi \in  B_{L^{2}}(0,2)$,
and
\begin{equation} \label{TFI-Der}
\frac{\partial F}{\partial \sigma} ( \sigma(\Psi),\Psi )
= 1 - 2 \theta'( \mathcal{V}_{\sigma(\Psi),N,\epsilon}(\Psi) )
\frac{\partial}{\partial \sigma}
\Big[ \mathcal{V}_{\sigma,N,\epsilon}(\Psi) \Big]_{\sigma(\psi)}
\geqslant \frac{1}{2}.
\end{equation}
Indeed, for $\sigma_{0} \in [0,\|\theta\|_{L^{\infty}}]$
and $\Psi \in B_{L^{2}}(0,2)$, we have
$$\frac{\partial}{\partial \sigma}
\Big[ \mathcal{V}_{\sigma,N,\epsilon}(\Psi) \Big]_{\sigma_{0}}
=
-2 \sum_{k=1}^{N} a_{k} \Re \left(
\langle \Psi , \frac{d\phi_{k,\sigma}}{d\sigma} \Big|_{\sigma_{0}} \rangle
\overline{\langle \Psi , \phi_{k,\sigma_{0}} \rangle}
\right)$$
where $a_{1}:=1$ and $a_{k}:=1-\epsilon$ for $k=2,..,N$.
Thus, using (\ref{vep-k-1}), we get
$$\Big| \frac{\partial}{\partial \sigma}
\Big[ \mathcal{V}_{\sigma,N,\epsilon}(\Psi) \Big]_{\sigma_{0}}
\Big| \leqslant 8NC^{*}.$$
We get the inequality in (\ref{TFI-Der})
thanks to the previous inequality and (\ref{hyp-theta}).

For every $\Psi \in B_{L^{2}}(0,2)$, the implicit function theorem
provides the existence of a local $C^{\infty}$ parameterization
$\tilde{\sigma}(\xi)$ for the solutions of $F(\sigma(\xi),\xi)=0$,
in a neighborhood of $\Psi$. The uniqueness of the fixed point
$\sigma(\xi)$ justifies that $\sigma$ and $\tilde{\sigma}$ coincide,
thus $\sigma$ is $C^{\infty}$. \endproof

\textbf{Acknowledgments :} The authors thank J-M. Coron, R. van
Handel, O. Kavian and P. Rouchon for helpful discussions.


\end{document}